\documentclass[11pt,letterpaper]{article}

\usepackage{url}
\usepackage{graphicx}
\usepackage{graphicx,caption,subcaption}
\usepackage{color}
\usepackage{epsfig}
\usepackage{latexsym}

\usepackage{amssymb}
\usepackage{amsmath}
\usepackage{amsthm}
\usepackage{framed}
\usepackage{a4wide}
\usepackage{enumerate}
\usepackage{makeidx}
\usepackage{bm}
\usepackage{epsfig}
\usepackage{graphics}
\usepackage{colortbl}
\usepackage{colordvi}
\usepackage{verbatim}
\usepackage{amsmath}
\usepackage{multirow} 
\usepackage{booktabs}

\newcommand{\probconverge}{\mbox{~$\stackrel{P}{\longrightarrow}$}}

\newcommand {\R} {\mathbb{R}}

\newcommand {\C} {{\rm I\kern-5.5pt C}}

\newcommand{\bP}[1]{{\mathbb{P}}\left[{#1}\right]}
\newcommand{\bE}[1]{{\mathbb{E}}\left[{#1}\right]}

\newcommand{\1}[1]{{\bf 1}\left[#1\right]}       

\newcommand{\fsquare}{\vrule height6pt width7pt depth1pt}   
\newcommand{\myproof}{{\hfill \\ \bf Proof. \ }}           
\newcommand{\myendpf}{\hfill\fsquare \\[0.1in]}             
\newcommand{\myvec}[1]{{\mbox{\boldmath{$#1$}}}}

\newtheorem{theorem}{Theorem}[section]

\newtheorem{lemma}[theorem]{Lemma}
\newtheorem{proposition}[theorem]{Proposition}
\newtheorem{corollary}[theorem]{Corollary}

\newtheorem{assumption}{Assumption}

\begin{document}

\sloppy
\date{}

\title{Asymptotic degree distributions\\
         in random threshold graphs
\thanks{This document does not contain technology or technical data controlled 
under either the U.S. International Traffic in Arms Regulations or the U.S. Export Administration Regulations. }
}
\author{Armand M. Makowski\thanks{A. M. Makowski is with the
Department of Electrical and Computer Engineering, and the
Institute for Systems Research, University of Maryland, College
Park, MD 20742 USA (e-mail: armand@isr.umd.edu).}
~and 
Siddharth Pal\thanks{S. Pal was with the Department of Electrical and Computer Engineering and 
and the Institute for Systems Research, University of Maryland, College
Park, MD 20742 USA. He is now with Raytheon BBN Technologies (email: siddharth.pal@raytheon.com). }
}
\maketitle

\noindent
\hrulefill \\
\normalsize We discuss several limiting degree distributions for a class of random threshold graphs in the many node regime. 
This analysis is carried out under a weak assumption on the distribution of the underlying fitness variable.
This assumption, which is satisfied by the exponential distribution, determines a natural scaling under which the following limiting results are shown:
The nodal degree distribution, i.e., the distribution of {\em any} node, converges in distribution to a limiting pmf. 
However, for each $d=0,1, \ldots $, the {\em fraction} of nodes with given degree $d$ converges only in distribution to a {\em non-degenerate}
random variable $\Pi(d)$ (whose distribution depends on $d$),
and {\em not} in probability to the aforementioned limiting nodal pmf as is customarily expected. 
The distribution of $\Pi(d)$  is identified only through its characteristic function.
Implications of this result include:
(i) The empirical node distribution may not be used as a proxy for or as an estimate to  the limiting nodal pmf;
(ii) Even in homogeneous graphs, 
the network-wide degree distribution and the nodal degree distribution may capture vastly different information;
and
(iii) Random threshold graphs with exponential distributed fitness do {\em not} provide 
an alternative scale-free model to the Barab\'asi-Albert model as was  argued by some authors;
the two models cannot be meaningfully compared in terms of their degree distributions!


\noindent
\hrulefill \\


\vfill

\pagebreak

\section{Introduction}
\label{sec:Introduction}

Graphs as network models are routinely studied through their degree distributions, and much of the attention
has focused on the {\em empirical} degree distribution that records the fractions of nodes with given degree value.
This distribution, which is easy to obtain from network measurements,
has been  found in many networks to obey a {\em power law} \cite[Section 1.4]{Durrett_Book}: If the network
comprises a large number $n$ of nodes and there are $N_n(d)$ nodes with degree $d$ among them, 
then the data reveals a behavior of the form
\begin{equation}
\frac{ N_n(d) }{n} \simeq C d^{-\alpha}
\label{eq:PowerLaw}
\end{equation}
for some $\alpha$ in the range $[2,3]$ (although there are occasional exceptions) and $C>0$
\cite{BarabasiAlbert, NewmanSurvey}.
See the monograph \cite[Section 4.2]{Durrett_Book} for an introductory discussion and references.
Statements such as (\ref{eq:PowerLaw})
are usually left somewhat vague as the range for $d$ is never carefully specified (in relation to $n$);
networks where (\ref{eq:PowerLaw}) was observed are often said to be {\em scale-free}.

The Barab\'asi-Albert model came to prominence 
as  the first random graph model to formally demonstrate 
the possibility of power law degree distribution in large networks \cite{BarabasiAlbert}:
The original Barab\'asi-Albert model is a growth model  which relies on the mechanism of {\em preferential attachment} --
Newly arriving nodes attach themselves to existing nodes with a probability proportional to their degrees at the time of arrival. 
As the number $n$ of nodes increases, Bollob\'as et al. \cite{Bollobas_BAmodel} proved that
\begin{equation}
           \frac{N_n(d)}{n} \probconverge_n ~ p_{\rm BA}(d), 
           \quad d=0,1, \ldots
           \label{eq:AsymptoticEmpiricalPMF_BA}
\end{equation}
where 
$\left (  \frac{N_n(d)}{n},\ d=0,1,\ldots \right )$ is the empirical degree distribution of the graph with $n$ nodes,
and the limiting pmf $\myvec{p}_{\rm BA} = \left ( p_{\rm BA} (d),\ d=0,1,\ldots \right )$ on $\mathbb{N}$ has the power-tail behavior
\begin{equation}
p_{\rm BA} (d) \sim d^{-3} 
\quad \mbox{($d \rightarrow \infty$)}.
\label{eq:PowerLaw1}
\end{equation}
Many generalizations of the Barab\'asi-Albert model have been proposed over the years:
Typically the convergence (\ref{eq:AsymptoticEmpiricalPMF_BA}) still holds for some limiting
pmf $\myvec{p} = \left ( p (d),\ d=0,1,\ldots \right )$ on $\mathbb{N}$
with (\ref{eq:PowerLaw1}) replaced by $p(d) \sim d^{-\tau} $ ($d \rightarrow \infty$)
for some $\tau > 0$. The various models distinguish themselves from each other  
by their ability to achieve a value $\tau$ in a particular range \cite[Section 4.2]{Durrett_Book}.

Although in some contexts preferential attachment is a reasonable assumption, 
it is predicated on the degree of existing nodes being available to newly arriving nodes. 
There are many situations where this assumption is questionable, 
and where the creation of a link between two nodes may instead result
in a mutual benefit based on their intrinsic attributes, e.g., authority, 
friendship, social success, wealth, etc.
{\em Random threshold graph} models, which were proposed by Caldarelli et al. \cite{CCDM},
incorporate this viewpoint in its simplest form as follows:
Let  $\{ \xi, \xi_k,\ k=1,2,\ldots  \}$ denote i.i.d. $\R_+$-valued random variables (rvs) with $\xi_k$ 
expressing the \lq\lq fitness" level associated with node $k$.
With $n$ nodes and a threshold $\theta > 0$, the random threshold graph $\mathbb{T}(n,\theta)$ 
postulates that two distinct nodes $i$ and $j$ form a connection (hence there is an undirected edge between them) if
\[
\xi_i +\xi_j > \theta,
\quad
\begin{array}{c}
i \neq j \\
i,j=1, \ldots , n. \\
\end{array}
\]

Interest in random threshold graphs has been spurred by the following observations:
The distribution of the degree rvs $D_{n,1}(\theta),\ldots, D_{n,n}(\theta)$ in $\mathbb{T}(n;\theta)$ 
(given by (\ref{eq:DegreeDefn})) is the same for all nodes. It is therefore appropriate 
to speak of {\em the} degree distribution of a node in $\mathbb{T}(n;\theta)$, namely that of $D_{n,1}(\theta) $.
Now consider the case when the fitness variable $\xi$ is  exponentially distributed with parameter $\lambda > 0$,
and  the threshold $\theta$ is scaled with the number $n$ of nodes according to 
the scaling $\theta^\star : \mathbb{N}_0 \rightarrow \mathbb{R}_+$ given by
\begin{equation}
\theta^\star_n = \lambda^{-1} \log n,
\quad n=2,3, \ldots 
\label{eq:ExponentialRTgraphs_scaling}
\end{equation}
In that setting  Fujihara et al. \cite{FIKMMU} have shown the distributional convergence
\begin{equation}
D_{n,1}(\theta^\star_n) \Longrightarrow_n D
\label{eq:DegreeConvergence_RTGraphsA}
\end{equation}
where the limiting rv $D$ has pmf $\myvec{p}_{\rm Fuj} = \left ( p_{\rm Fuj} (d),\ d=0,1,\ldots \right )$ on $\mathbb{N}$ 
with power-tail behavior
\begin{equation}
p_{\rm Fuj}(d) \sim d^{-2}
\quad \mbox{($d \rightarrow \infty$)}.
\label{eq:PowerLawFujihara}
\end{equation}
The result (\ref{eq:DegreeConvergence_RTGraphsA})-(\ref{eq:PowerLawFujihara})
has led some researchers \cite{CCDM,SC} to conclude that random threshold graphs
can model scale-free networks (albeit with $\tau =2$) without having to resort to either a growth process or a preferential attachment mechanism,
and as such they provide  an {\em alternative} to the Barab\'asi-Albert model.
However, a moment of reflection should lead one to question
this conclusion given the evidence available so far. 
Indeed, the statement (\ref{eq:AsymptoticEmpiricalPMF_BA})
concerns an empirical  degree distribution which is computed network-wide, whereas the convergence
(\ref{eq:DegreeConvergence_RTGraphsA})-(\ref{eq:PowerLawFujihara})
addresses the distributional behavior of the degree of a {\em single} node, 
its distribution being identical across nodes.

A natural question is whether this discrepancy can be resolved in the large network limit.
More precisely, for each $d=0,1, \ldots $, let $N_n(d;\theta)$ denote the number of nodes in $\mathbb{T}(n;\theta)$ which have degree $d$, namely
\[
N_n(d;\theta) = \sum_{k=1}^n \1{ D_{n,k}(\theta)= d },
\quad n=2,3, \ldots 
\]
In analogy with (\ref{eq:AsymptoticEmpiricalPMF_BA}), 
is it indeed the case that 
\begin{equation}
 \frac{ N_n(d;\theta^\star_n) }{ n }  \probconverge_n ~  p_{\rm Fuj}(d)
 \label{eq:Intro_ConvergenceProbability}
\end{equation}
where the pmf $\myvec{p}_{\rm Fuj}$ is the one appearing at (\ref{eq:DegreeConvergence_RTGraphsA})-(\ref{eq:PowerLawFujihara})?
Only then would random threshold graphs (under exponentially distributed fitness) 
be confirmed as a {\em bona fide} scale-free alternative model to the Barab\'asi-Albert model (as described by
(\ref{eq:AsymptoticEmpiricalPMF_BA})-(\ref{eq:PowerLaw1})).

In this paper, for each $d=0,1, \ldots $, we show  that there exists a {\em non-degenerate} $[0,1]$-valued rv $\Pi(d)$ such that
\begin{equation}
 \frac{ N_n(d;\theta^\star_n) }{ n } \Longrightarrow_n \Pi(d)
 \label{eq:Intro_Weakconvergence}
\end{equation}
where the scaling $\theta^\star : \mathbb{N}_0 \rightarrow \mathbb{R}_+$ is the one defined at (\ref{eq:ExponentialRTgraphs_scaling}) -- In fact
we establish such a result for a very large class of fitness distributions (with the scaling $\theta^\star : \mathbb{N}_0 \rightarrow \mathbb{R}_+$ 
modified accordingly).
The non-degeneracy of the rv $\Pi(d)$ in (\ref{eq:Intro_Weakconvergence}) implies 
that  (\ref{eq:Intro_ConvergenceProbability}) {\em cannot}  hold, and
random threshold graphs with exponential distributed fitness do {\em not} provide an alternative scale-free model to the Barab\'asi-Albert model
(as understood by (\ref{eq:AsymptoticEmpiricalPMF_BA})).
Only the convergence (\ref{eq:AsymptoticEmpiricalPMF_BA}) has meaning  in the preferential attachment model
while  the convergence (\ref{eq:DegreeConvergence_RTGraphsA}) has no equivalent there, 
the situation being reversed for random threshold graphs -- 
The two models cannot be meaningfully compared in terms of their degree distributions!
Thus, {\em even} in homogeneous graphs, 
the network-wide degree distribution and the nodal degree distribution may capture vastly different information.  
This issue was also investigated more broadly by the authors in the references
\cite{SPal_Thesis, PalMakowski_CDC2015, PalMakowski-TNSE}; see comments following Corollary \ref{Cor:MainResult}.

We close with a summary of the contents of the paper:
Random threshold graphs are introduced in Section \ref{sec:RTGs} together with the needed notation
and assumptions.
As we consider situations that generalize the case of exponentially distributed fitness rvs,
the scaling (\ref{eq:ExponentialRTgraphs_scaling})  is now replaced 
by a scaling $\theta^\star : \mathbb{N}_0 \rightarrow \mathbb{R}_+$ satisfying Assumption \ref{ass:A}.
This assumption is determined by the probability distribution of $\xi$, and ensures the convergence
$D_{n,1}(\theta^\star_n) \Longrightarrow_n D$ for some limiting rv $D$ which is conditionally Poisson (given $\xi$)
[Proposition \ref{prop:OneDimConvergence}].
Section \ref{sec:MainResults} presents the main result of the paper [Theorem \ref{thm:MainResult}],
namely that in the setting of Section \ref{sec:RTGs}, under Assumption \ref{ass:A},
the distributional convergence (\ref{eq:Intro_Weakconvergence}) holds with a non-degenerate limit identified only through its characteristic function
(\ref{eq:Phi_d(t)_Defn})-(\ref{eq:phi_d(t)}).
A proof of Theorem \ref{thm:MainResult} is given in  Section \ref{sec:ProofMainResult} and is rooted in the method of moments
[via Proposition \ref{prop:ConvergenceMoments} established in Section \ref{sec:ProofPropConvergenceMoments}].
The main technical step is contained in Proposition \ref{prop:MultiDimConvergence}; 
the proof of this multi-dimensional version of Proposition \ref{prop:OneDimConvergence} is given 
in several steps which are presented from Section \ref{sec:ProofMultDimConvergence} to Section \ref{sec:ProofPropLimitPGFsPart_III}.
Section \ref{sec:Simulations} illustrates through limited simulations the failure of  (\ref{eq:Intro_ConvergenceProbability})
and the validity of  (\ref{eq:Intro_Weakconvergence}) in the case of  random threshold graphs with exponentially distributed fitness.

\section{Random threshold graphs}
\label{sec:RTGs}

First some notation and conventions:
The random variables (rvs) under consideration are all 
defined on the same probability triple $(\Omega, {\cal F}, \mathbb{P})$. 
The construction of a sufficiently large probability triple carrying all needed rvs 
is standard and omitted in the interest of brevity.
All probabilistic statements are made with respect to the probability
measure $\mathbb{P}$, and we denote the corresponding expectation operator by $\mathbb{E}$. 
The notation $\probconverge_n$ (resp.  $\Longrightarrow_n$)  is used to signify convergence in probability
(resp. convergence in distribution) (under $\mathbb{P}$) with $n$ going to infinity;
see the monographs \cite{BillingsleyBook, ChungBook, Shiryayev}
for definitions and properties.
If $E$ is a subset of $\Omega$, then $\1{E}$ denotes the indicator
of the set $E$ with the usual understanding that $\1{E}(\omega) = 1$ (resp. $\1{E}(\omega) = 0$) if $\omega \in E$
(resp. $\omega \notin E$). 
The symbol $\mathbb{N}$ (resp. $\mathbb{N}_0$)  denotes the set
of non-negative (resp. positive) integers.

\subsection{Model}
\label{subsec:RTGs+Model}

The setting is that of \cite{MakowskiYagan-JSAC}:
Let $\{ \xi, \xi_k, \ k=1,2, \ldots \}$
denote a collection of i.i.d. $\mathbb{R}_+$-valued rvs defined on 
the probability triple $(\Omega, {\cal F}, \mathbb{P})$,
each distributed according to a given
(probability) distribution function 
$F: \mathbb{R} \rightarrow [0,1]$.
With $\xi$ acting as  a generic representative for this sequence of i.i.d. rvs,
we have
\[
\bP{ \xi \leq x } = F(x),
\quad x \in \mathbb{R}.
\]
At minimum we assume that $F$ is a {\em continuous} function on $\mathbb{R}$ with support on $[0,\infty)$, namely
\begin{equation}
F(x) = 0, \quad x \leq 0.
\label{eq:NonNegative}
\end{equation}

Once $F$ is specified, random thresholds graphs are characterized by two parameters, 
namely the number $n$ of nodes and a threshold value $\theta > 0$:
The network comprises $n$ nodes, labelled $k=1, \ldots , n$,
and to each node $k$ we assign a {\em fitness} variable (or weight) $\xi_k$
For distinct $k,\ell =1, \ldots , n$,
nodes $k$ and $\ell$ are declared to be adjacent if
\begin{equation}
\xi_k + \xi_\ell > \theta ,
\label{eq:AdjacencyDefnRTGs}
\end{equation} 
in which case we say that an undirected link exists between these two nodes.
The {\em random threshold} graph is the (undirected) random graph $\mathbb{T}(n;\theta)$
on the set of vertices $\{ 1, \ldots , n \}$
defined by the adjacency notion (\ref{eq:AdjacencyDefnRTGs}).

For each $k=1,2, \ldots , n$,
the degree of node $k$ in $\mathbb{T}(n;\theta)$
is the rv $D_{n,k}(\theta)$ given by
\begin{equation}
D_{n,k}(\theta)
= \sum_{\ell=1, \ \ell \neq k}^n \1{ \xi_k + \xi_\ell > \theta } .
\label{eq:DegreeDefn}
\end{equation}
Under the enforced independence assumptions, conditionally on $\xi_k$, the rv $D_{n,k}(\theta)$ is a Binomial rv 
${\rm Bin}(n-1; 1-F(\theta-\xi_k) )$.
The rvs $D_{n,1}(\theta), \ldots , D_{n,n}(\theta)$ being exchangeable,
let $D_n(\theta)$ denote any $\mathbb{N}$-valued rv which is distributed according to
their common pmf. 

\subsection{Existence of a limiting degree distribution}
\label{subsec:RTGs+ExistenceLimitDegree}

Throughout we make the following assumption on $F$.

\begin{assumption}
{\sl 
There exists a scaling $\theta^\star: \mathbb{N}_0 \rightarrow \mathbb{R}_+$ with the property
\begin{equation}
\lim_{n \rightarrow \infty } \theta^\star_n = \infty ,
\label{eq:GoingToInfinity}
\end{equation}
such that
\begin{equation}
\lim_{n \rightarrow \infty} n \left ( 1 - F(\theta^\star_n - x ) \right ) = \lambda(x),
\quad x \geq 0
\label{eq: KeyCondition1}
\end{equation}
for some non-identically zero mapping $\lambda : \mathbb{R}_+ \rightarrow \mathbb{R}_+$.
}
\label{ass:A}
\end{assumption}

The mapping $\lambda : \mathbb{R}_+ \rightarrow \mathbb{R}_+$ is necessarily non-decreasing.
The following result overlaps with a similar result by
Fujihara et al. \cite[Thm. 2, p. 362]{FIKMMU}.

\begin{proposition}
{\sl Under Assumption \ref{ass:A}, there exists an $\mathbb{N}$-valued rv $D$ such that
\begin{equation}
D_n(\theta^\star_n) \Longrightarrow_n D.
\label{eq:OneDimConvergence}
\end{equation}
The rv $D$ is conditionally Poisson with pmf given by
\begin{equation}
\bP{ D = d } = \bE{ \frac{ \lambda (\xi)^d}{d!} e^{-\lambda (\xi) } },
\quad d=0,1, \ldots
\label{eq:ConditionalPoisson}
\end{equation}
}
\label{prop:OneDimConvergence}
\end{proposition}
The convergence (\ref{eq:OneDimConvergence}) is equivalent to
\begin{equation}
\lim_{n \rightarrow \infty} \bP{ D_n(\theta^\star_n) = d } = \bP{ D = d },
\quad d=0,1, \ldots
\label{eq:OneDimConvergence2}
\end{equation}
If the mapping $\lambda : \mathbb{R}_+ \rightarrow \mathbb{R}_+$ assumes a constant value $c>0$, i..e.,
$\lambda(x) = c $ for all $x \geq 0$, then the  rv $D$ is a Poisson rv with parameter $c$.

\myproof
Fix $n=2,3, \ldots $, $\theta > 0$ and $z$ in $\mathbb{R}$.
Standard pre-conditioning arguments yield
\begin{eqnarray}
\bE{z^{D_{n}(\theta) } }
=
\bE{ \prod_{\ell=2}^n z^{ \1{ \xi_1 + \xi_\ell > \theta } } }
=
\bE{ 
\left ( 
\bE{ z^{ \1{ x + \xi > \theta } }  }_{x = \xi_1} 
\right)^{n-1} }
\nonumber 
\end{eqnarray}
under the enforced independence assumptions where
\begin{eqnarray}
\bE{ z^{ \1{ x + \xi > \theta } }  }
&=&
\bP{ x + \xi \leq \theta }  + \bP{  x + \xi > \theta }  z
\nonumber \\
&=&
1 - (1-z) \bP{  x + \xi > \theta }
\nonumber \\
&=&
1 - (1-z)  \left ( 1 - F( \theta - x ) \right ),
\quad x \in \mathbb{R}.
\end{eqnarray}
Consequently, upon using (\ref{eq:NonNegative}), we get
\begin{eqnarray}
\bE{z^{D_{n}(\theta) } }
&=&
\bE{
\left ( 1 - (1-z)  \left ( 1 - F( \theta - \xi ) \right ) \right )^{n-1}
}
\nonumber \\
&=&
z^{n-1}  \bP{ \xi > \theta } 
+
\bE{ \1{ \xi \leq \theta } \left ( 1 - (1-z)  \left ( 1 - F( \theta - \xi ) \right ) \right )^{n-1} }.
\label{eq:Equality1}
\end{eqnarray}

Now replace $\theta$ by $\theta^\star_n$ in (\ref{eq:Equality1})
according to the scaling $\theta^\star: \mathbb{N}_0 \rightarrow \mathbb{R}_+$ stipulated in Assumption \ref{ass:A}, and 
let $n$ go to infinity in the resulting equality when $|z| \leq 1$: 
It is plain that $\lim_{n \rightarrow \infty} \bP{ \xi > \theta^\star_n } z^{n-1} =0$, while standard arguments \cite[Prop. 3.1.1., p. 116]{EKM}
yield
\[
\lim_{n \rightarrow \infty} 
\1{ \xi \leq \theta^\star_n} \left ( 1 - (1-z)  \left ( 1 - F( \theta^\star_n - \xi ) \right ) \right )^{n-1} 
= e^{ -(1-z) \lambda(\xi) }
\]
since
$\lim_{n \rightarrow \infty }  n (1-z)  \left ( 1 - F( \theta^\star_n - \xi ) \right ) = (1-z) \lambda(\xi)$
under Assumption \ref{ass:A}.
Invoking  the Bounded Convergence Theorem we obtain
\[
\lim_{n \rightarrow \infty} \bE{ z^{D_n(\theta^\star_n)} }
= \bE{ e^{-(1-z) \lambda(\xi) }},
\quad |z| \leq 1
\]
and the desired conclusion 
(\ref{eq:OneDimConvergence})--(\ref{eq:ConditionalPoisson})
follows by standard arguments upon noting that the right-hand side is the probability generating function (pgf) of the pmf (\ref{eq:ConditionalPoisson}).
\myendpf

Assumption  \ref{ass:A} holds in a number of interesting cases;
in what follows we use the standard notation $x^+ = \max(x,0)$ for $x$ in $\mathbb{R}$:
When $\xi$ is exponentially distributed with parameter $\lambda > 0$, 
namely
\begin{equation}
\bP{ \xi > x } =e^{-\lambda x^+ },
\quad x \in \mathbb{R},
\label{eq:ExponentialFitness}
\end{equation}
Assumption  \ref{ass:A} holds with
\begin{equation}
\lambda(x) = e^{ \lambda x },
\quad x \geq 0
\label{eq:LambdaExponential}
\end{equation}
if we take  $\theta^\star_n = \lambda^{-1} \log n$ for all $n =1,2, \ldots $.
In this case, the pmf of $D$ is the pmf $\myvec{p}_{\rm Fuj}$ 
appearing at (\ref{eq:DegreeConvergence_RTGraphsA})-(\ref{eq:PowerLawFujihara}); it is given by
\begin{equation}
p_{\rm Fuj} (d)  = \bP{ D = d } = \bE{ \frac{e^{d \lambda \xi} }{d!} e^{-e^{\lambda \xi} } } 
= \frac{1}{d!} \int_0^\infty e^{dx} e^{-e^{x} } e^{-x} dx
\quad d=0,1, \ldots
\label{eq:ConditionalPoisson2}
\end{equation}
as we substitute (\ref{eq:LambdaExponential}) into the expression (\ref{eq:ConditionalPoisson}).
Note that $p_{\rm Fuj} (d)$ does not depend on $\lambda$ since $\lambda \xi$ is exponentially distributed with unit parameter
if $\xi$ is exponentially distributed with parameter $\lambda$.

The second case deals with heavy-tailed rvs:
The rv $\xi$ is said to be a Pareto rv with parameters $\nu > 0$ and $a > 0$ if
\begin{equation}
\bP{ \xi > x } = \left ( \frac{a}{a+x^+} \right )^\nu,
\quad x \in \mathbb{R}.
\label{eq:ParetoFitness}
\end{equation}
Assumption  \ref{ass:A} holds with
$\theta^\star_n = an^{\frac{1}{\nu}}$ for all $n =1,2, \ldots$, and $ \lambda(x) = 1 $ for all $x \geq 0$, in which case $D$ is a Poisson rv with unit parameter.

\section{Main results}
\label{sec:MainResults}

Fix $n=2,3, \ldots $ and $\theta > 0$.
For each $d=0,1, \ldots $, the rv $N_n(d;\theta)$ defined by
\begin{equation}
N_n(d;\theta) = \sum_{k=1}^n \1{ D_{n,k}(\theta)= d }
\label{eq:CountVariable}
\end{equation}
counts the number of nodes in $\{1, \ldots , n \}$ which have degree $d$ in $\mathbb{T}(n; \theta)$. 
The fraction of nodes in $\{1, \ldots , n \}$ with degree $d$ in $\mathbb{T}(n;\theta)$ is then given by
\begin{equation}
P_n (d;\theta) = \frac{ N_n(d;\theta) }{n}.
\label{eq:EmpiricalPMF}
\end{equation}
The main result of the paper is concerned with the following convergence.

\begin{theorem}
{\sl Assume Assumption \ref{ass:A} to hold.
For each $d=0,1, \ldots $, there exists a non-degenerate $[0,1]$-valued rv $\Pi(d)$ such that
\begin{equation}
P_n (d;\theta^\star_n)  \Longrightarrow_n \Pi(d)
 \label{eq:MainResultConvergence}
\end{equation}
where the scaling $\theta^\star : \mathbb{N}_0 \rightarrow \mathbb{R}_+$ is the one postulated in
Assumption \ref{ass:A}. Furthermore, it holds that $\bE{ \Pi(d) } = \bP{ D=d } $ and ${\rm Var} [ \Pi (d) ] > 0$
with $D$ being the limiting rv whose existence is established in Proposition \ref{prop:OneDimConvergence}.
}
\label{thm:MainResult}
\end{theorem}

In the course of proving Theorem \ref{thm:MainResult} (in Section \ref{sec:ProofMainResult}), we determine
the distribution of the rv $\Pi(d)$ through its characteristic function (\ref{eq:phi_d(t)}).
The non-degeneracy of the rv $\Pi(d)$ has the following consequence.

\begin{corollary}
{\sl Assume Assumption \ref{ass:A} to hold.
For each $d=0,1, \ldots $, the sequence $\left \{  P_n(d;\theta^\star_n) , \ n=1,2, \ldots \right  \}$
cannot converge in probability to a constant, i.e., there exists no constant $L(d)$ such that
\begin{equation}
P_n(d;\theta^\star_n)  \probconverge_n ~ L(d).
\label{eq:MainResultNonConvergence}
\end{equation}
}
\label{Cor:MainResult}
\end{corollary}
Corollary \ref{Cor:MainResult} was announced in the conference paper \cite{PalMakowski_CDC2015}
when the fitness variables are exponentially distributed; in \cite{PalMakowski-TNSE}  
the failure of the convergence (\ref{eq:MainResultNonConvergence}) was shown in the exponential  case
with the help of asymptotic properties of order statistics.
Here, a fuller picture is obtained:
 Corollary \ref{Cor:MainResult} is a by-product of the weak convergence  (\ref{eq:MainResultConvergence})
(which replaces  the non-convergence (\ref{eq:MainResultNonConvergence}) and requires only Assumption  \ref{ass:A} to hold),
and of the non-degenerate nature of the  limiting rv $\Pi (d)$.

The remainder of the paper is concerned with establishing Theorem \ref{thm:MainResult};
its proof relies on  the method of moments, and proceeds through
Proposition \ref{prop:MultiDimConvergence} and Proposition \ref{prop:ConvergenceMoments}
which are stated below.
In Section \ref{sec:ProofMainResult}  we rely on these two intermediary results to construct a short proof of Theorem \ref{thm:MainResult}.
Proposition \ref{prop:MultiDimConvergence} 
contains a  multi-dimensional version of Proposition \ref{prop:OneDimConvergence}, and
provides the core technical content behind Theorem \ref{thm:MainResult};
a multi-step proof is presented  from Section \ref{sec:ProofMultDimConvergence} to Section \ref{sec:ProofPropLimitPGFsPart_III}.

\begin{proposition}
{\sl
Assume Assumption \ref{ass:A} to hold.
For each $r=1, 2, \ldots$, 
there exists an $\mathbb{N}^r$-valued rv $(D_1, \ldots ,D_r)$
such that
\begin{equation}
\left ( D_{n,1}(\theta^\star_n) , \ldots , D_{n,r}(\theta^\star_n) \right )
\Longrightarrow_n (D_1, \ldots ,D_r).
\label{eq:MultiDimConvergence=r}
\end{equation}
The limiting rvs $D_1, \ldots , D_r$ are exchangeable, but not independent,
each being distributed according to the
limiting rv $D$ whose existence is established in Proposition \ref{prop:OneDimConvergence}.
}
\label{prop:MultiDimConvergence}
\end{proposition}

The next step,  established in Section \ref{sec:ProofPropConvergenceMoments},
deals with the needed convergence to apply the method of moments.

\begin{proposition}
{\sl
Assume Assumption \ref{ass:A} to hold.
For each $r=1, 2, \ldots$,  we have
\begin{equation}
\lim_{n \rightarrow \infty}
\bE{ P_n(d;\theta^\star_n) ^r }
= \bP{ D_1 = d , \ldots , D_r = d },
\quad d=0,1, \ldots 
\label{eq:ConvergenceMoments=r}
\end{equation}
where the $\mathbb{N}^r$-valued rv $(D_1, \ldots ,D_r)$ is the limiting rv  whose existence was established in
Proposition \ref{prop:MultiDimConvergence}.
}
\label{prop:ConvergenceMoments}
\end{proposition}

\section{Simulation results}
\label{sec:Simulations}

In order to illustrate the difference between the convergence statements 
\eqref{eq:AsymptoticEmpiricalPMF_BA} and \eqref{eq:Intro_Weakconvergence},
we have carried out a limited
set of simulation experiments which are discussed in this section.
Throughout, the fitness variable $\xi$ is taken to be exponentially distributed with parameter $\lambda=1$, and the threshold is scaled 
in accordance with (\ref{eq:ExponentialRTgraphs_scaling}), namely
$\theta _n ^\star = \log n$ for each $n=2,3,\ldots $.

With the number $n$ of nodes given, we generate $R$ mutually independent versions of the random threshold graph $\mathbb{T}(n; \theta ^\star _n)$;
these realizations are denoted $ \mathbb{T}^{(1)}(n; \theta ^\star _n),  \mathbb{T}^{(2)}(n; \theta ^\star _n) ,  \ldots ,  \mathbb{T}^{(R)}(n; \theta ^\star _n)$.
For each $k=1,2,\ldots,n$ and $r=1,2,\ldots,R$, let $D^{(r)}_{n,k}(\theta^\star_n)$ denote the degree of node $k$ in
the random graph $\mathbb{T}^{(r)}(n; \theta ^\star _n)$, and for $d=0,1, \ldots $, let $N^{(r)}_n (d;\theta^\star_n)$ denote the 
number of nodes with degree $d$ in $\mathbb{T}^{(r)}(n; \theta ^\star _n)$. 

\begin{figure}[h]
\centering
\begin{subfigure}{.55\textwidth}
 \centering
 \includegraphics[width=.8\linewidth]{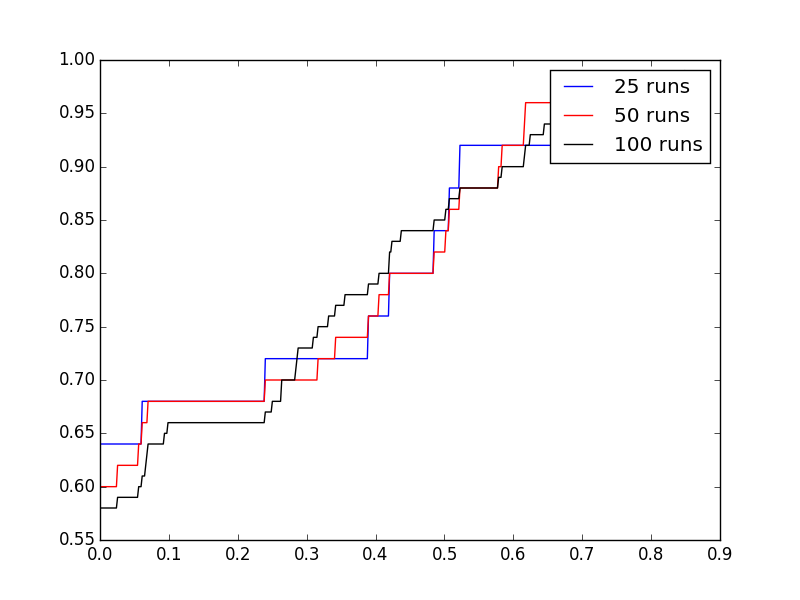}
 \caption{Varying $R$ and $n=30000$}
 \label{fig:RT_VaryRuns_deg0}
\end{subfigure}%
\begin{subfigure}{.55\textwidth}
 \centering
 \includegraphics[width=.8\linewidth]{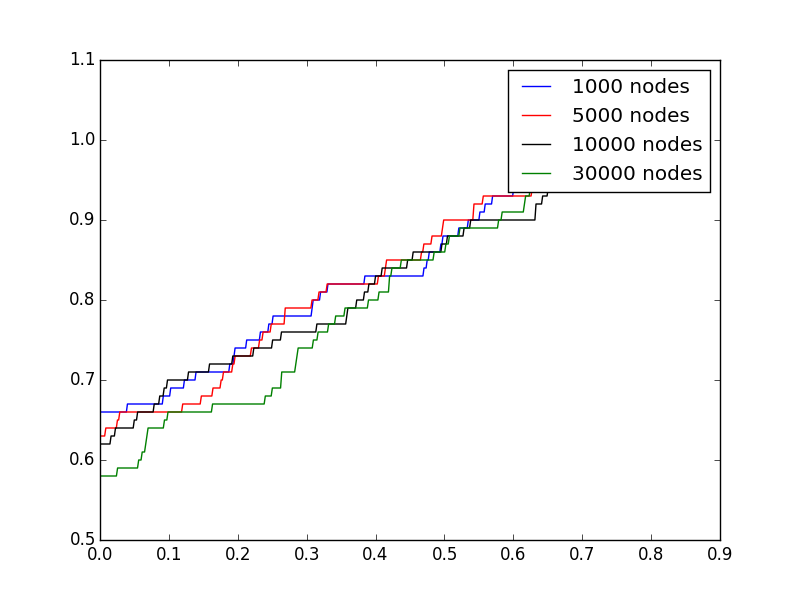}
 \caption{Varying $n$ and $R=100$}
 \label{fig:RT_VarySize_deg0}
\end{subfigure}
\caption{Histogram $H_{n,R}(d; .)$ for degree $d=0$ with varying number of nodes $n$ and the number of runs $R$ held fixed, and vice versa}
\label{fig:RT_deg0}
\end{figure}

\begin{figure}[h]
\centering
\begin{subfigure}{.55\textwidth}
 \centering
 \includegraphics[width=.8\linewidth]{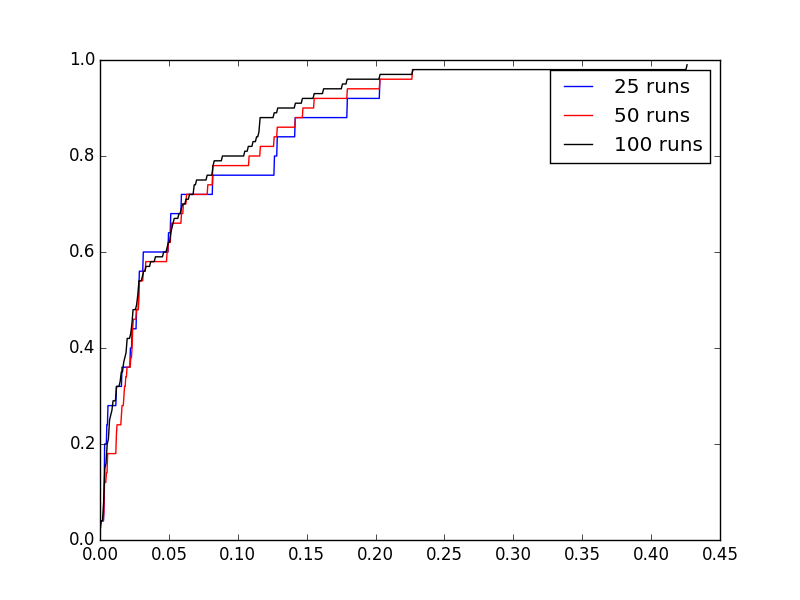}
 \caption{Varying $R$ and $n=30000$}
 \label{fig:RT_VaryRuns_deg5}
\end{subfigure}%
\begin{subfigure}{.55\textwidth}
 \centering
 \includegraphics[width=.8\linewidth]{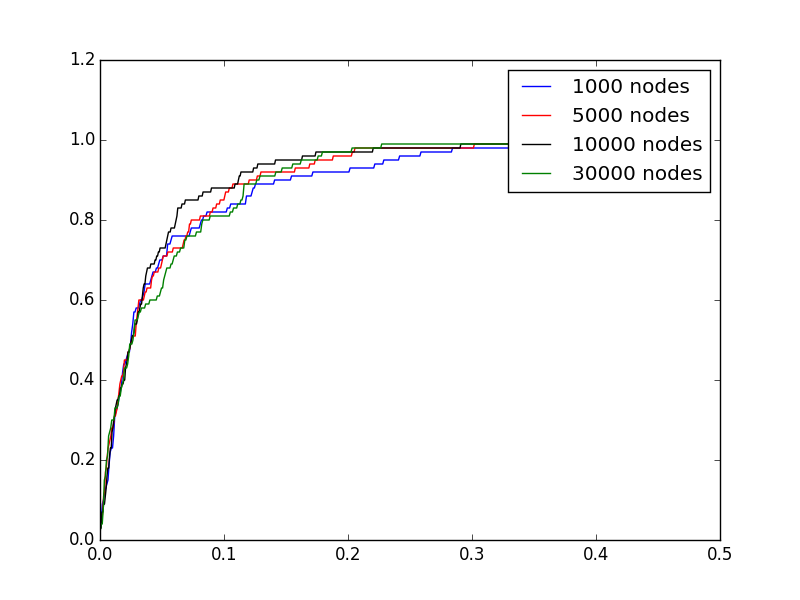}
 \caption{Varying $n$ and $R=100$}
 \label{fig:RT_VarySize_deg5}
\end{subfigure}
\caption{Histogram $H_{n,R}(d; .)$ for $d=5$ with varying number of nodes $n$ and the number of runs $R$ held fixed, and vice versa}
\label{fig:RT_deg5}
\end{figure}

\begin{figure}[h]
\centering
\begin{subfigure}{.55\textwidth}
 \centering
 \includegraphics[width=.8\linewidth]{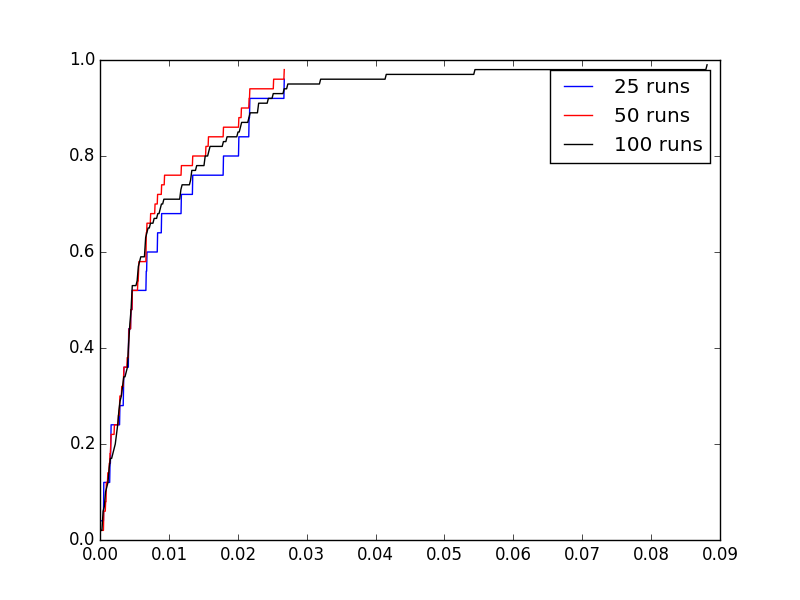}
 \caption{Varying $R$ and $n=30000$}
 \label{fig:RT_VaryRuns_deg10}
\end{subfigure}%
\begin{subfigure}{.55\textwidth}
 \centering
 \includegraphics[width=.8\linewidth]{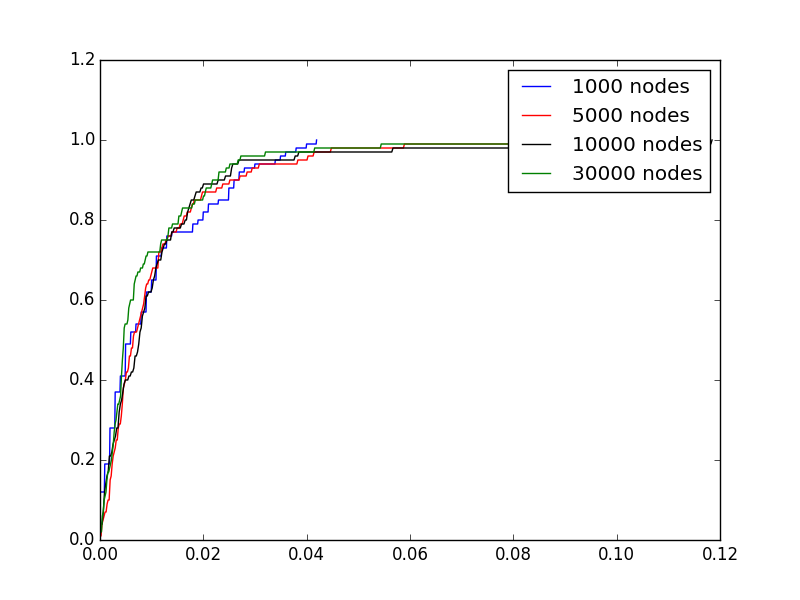}
 \caption{Varying $n$ and $R=100$}
 \label{fig:RT_VarySize_deg10}
\end{subfigure}
\caption{Histogram $H_{n,R}(d; .)$ for $d=10$ with varying number of nodes $n$ and the number of runs $R$ held fixed, and vice versa}
\label{fig:RT_deg10}
\end{figure}

\paragraph{The rv $\Pi(d)$ is non-degenerate}
Fix $d=0,1, \ldots $. 
On the strength of Theorem \ref{thm:MainResult}, a natural way to produce a reasonably good estimate for the probability distribution
of the rv $\Pi(d)$ is to follow a simple two-step procedure:
On the basis of the $R$ i.i.d. realizations of the random threshold graph $\mathbb{T}(n;\theta^\star_n)$,
a standard estimate of the probability distribution of $P_n(d;\theta_n ^\star)$ is provided by the  histogram
\[
H_{n,R} (d;x)  
=  \frac{1}{R}
\sum _{r=1} ^R \1{ \frac{N _n ^{(r)} (d;\theta _n ^\star)}{n} \leq x }, \ x \in \mathbb{R}.
\]
The probability distribution $x\rightarrow H_{n,R}(d; x)$ has support on $[0,1]$ with
$H_{n,R} (d;x)  = 0 $ for $x < 0$ and $H_{n,R} (d;x)  = 1 $ for $1 \leq x$, as does the probability distribution of $\Pi(d)$.
Under the enforced independence assumptions,  the Glivenko-Cantelli Theorem \cite[p. 103]{BillingsleyBook} asserts that
\begin{equation}
\lim_{ R \rightarrow \infty} 
\left (
\sup_{ 0 \leq x \leq 1 }
\left |
H_{n,R} (d;x) - \bP{ P_n(d;\theta_n ^\star) \leq x } 
\right |
\right )
= 0
\quad \mbox{a.s.}
\label{eq:Simul+0}
\end{equation}
Thus, for large $R$ (possibly dependent on $n$), the probability distribution 
of the rv $P_n(d;\theta_n ^\star)$ is uniformly well approximated by the histogram
$x\rightarrow H_{n,R} (d; x)$ with high probbability.

On the other hand, Theorem \ref{thm:MainResult} states that
\begin{equation}
\lim_{n \rightarrow \infty} \bP{ P_n(d;\theta_n ^\star) \leq x }  = \bP{ \Pi(d) \leq x },
\quad x \in \mathcal{C}( \Pi(d) )
\label{eq:Simul+2}
\end{equation}
where $\mathcal{C}( \Pi(d) )$ is the set of points of continuity of the probability distribution of $\Pi(d)$.
Thus, for each $x$ in $\mathcal{C}( \Pi(d) )$, the probability $\bP{ \Pi(d) \leq x }$ will be well approximated by $\bP{ P_n(d;\theta_n ^\star) \leq x } $
when $n$ is large (possibly dependent on $x$).

Combining (\ref{eq:Simul+0}) and (\ref{eq:Simul+2}) with a simple triangle inequality argument
naturally leads us to propose the approximation
\begin{equation}
\bP{ \Pi(d) \leq x } =_{\rm Approx}  H_{n,R}  (d;x),
\quad x \in \mathcal{C}( \Pi(d) )
\label{eq:Simul+5}
\end{equation}
with integers $n$ and $R$ selected sufficiently large.
Put differently, we expect the probability distribution of $\Pi(d)$ to be well approximated by the histogram
$x \rightarrow H_{n,R}  (d;x)$ if we select both $n$ and $R$ to be large. 
If such a histogram were found to be very different from a step function,
this would provide compelling evidence that  (\ref{eq:MainResultNonConvergence}) cannot hold, and that the rv $\Pi(d)$ is {\em not} degenerate.

In Figures~\ref{fig:RT_deg0}-\ref{fig:RT_deg10}, we show the approximating histogram $H_{n,R}(d;.)$ for the values $d=0,5,10$, with
a varying number $R$ of runs and a varying number $n$ of graph sizes.
Figure~\ref{fig:RT_deg0} deals with $d=0$: Figure~\ref{fig:RT_VaryRuns_deg0} shows the histograms for  $n=30000$ 
with increasing values $R=25,50,100$; the shape of the corresponding histograms do not change significantly. 
In Figure~\ref{fig:RT_VarySize_deg0}, with $R=100$, increasing the graph size 
$n=1000, 5000,10000, 30000$ also does not change the histograms significantly. This points to the non-degeneracy of $\Pi(0)$ 
since in all cases the approximating histograms are reasonably close together but never approximate, even remotely, a step function.
Figures~\ref{fig:RT_deg5} and \ref{fig:RT_deg10} exhibit histogram plots for $d=5,10$ under similar conditions;
the conclusions are identical to the ones reached in the case $d=0$, with the evidence being possibly even stronger since the histograms
appear to \lq\lq bend" in a concave manner.

\paragraph{Empirical degree distribution vs. nodal degree distribution}
As noted earlier, in the exponential case, the limiting rv $D$ appearing at (\ref{eq:OneDimConvergence})
has pmf $\myvec{p}_{\rm Fuj}$ given by  (\ref{eq:ConditionalPoisson2}), namely
\begin{equation}
p_{\rm Fuj} (d)  = \bP{ D = d } = \frac{1}{d!}  \int_0^\infty  \frac{e^{d x} }{d!} e^{-e^{x} } e^{-x} dx,
\quad d=0,1, \ldots
\label{eq:ConditionalPoisson3}
\end{equation}
For $d=0$ we numerically evaluate the appropriate integral with
\[
p_{\rm Fuj} (0)  = \int_0^\infty  e^{-e^{x} } e^{-x} dx = \int_0^1 e^{-t^{-1}} dt \simeq 0.1485.
\]
For $d = 2, 3, \ldots $ we note from (\ref{eq:ConditionalPoisson3}) that
\begin{eqnarray}
p_{\rm Fuj} (d)  &=&  \frac{1}{d!} \int_0^\infty  e^{(d-1) x} e^{-e^{x} } dx
\nonumber \\
&=&  \frac{1}{d!} \int_1^\infty  t^{d-2} e^{-t} dt
\quad \mbox{[Change of variable $t=e^x$]}
\nonumber \\
&=&  \frac{1}{d!} 
\left ( \int_0^\infty  t^{d-2} e^{-t} dt - \int_0^1  t^{d-2} e^{-t} dt \right ) 
\nonumber \\
&=&  \frac{1}{d!} 
\left ((d-2)!  - \int_0^1  t^{d-2} e^{-t} dt \right ) 
\nonumber
\end{eqnarray}
whence
\[
p_{\rm Fuj} (d) -  \frac{1}{d(d-1)} 
=
\frac{ \int_0^1  t^{d-2} e^{-t} dt }{d!} .
\]
The bound
\[
\varepsilon (d)
= \left | p_{\rm Fuj} (d)   -  \frac{1}{d(d-1)} \right | \leq  \frac{1}{d!}
\]
follows and readily yields (\ref{eq:PowerLawFujihara}).
This suggests the approximation 
\[
 p_{\rm Fuj} (d)   =_{\rm Approx} \frac{1}{d(d-1)}
\]
whose accuracy dramatically increases with $d$ increasing as $\varepsilon(d)$ decreases very rapidly, e.g.,
the approximations $ p_{\rm Fuj} (5)    =_{\rm Approx} \frac{1}{20}$ (with an error less than $1/120$)
and $ p_{\rm Fuj} (10)   =_{\rm Approx} \frac{1}{90}$ (with an error less than $1/10!$) are already tight.

\begin{figure}[h]
\centering
\begin{subfigure}{.55\textwidth}
 \centering
 \includegraphics[width=.8\linewidth]{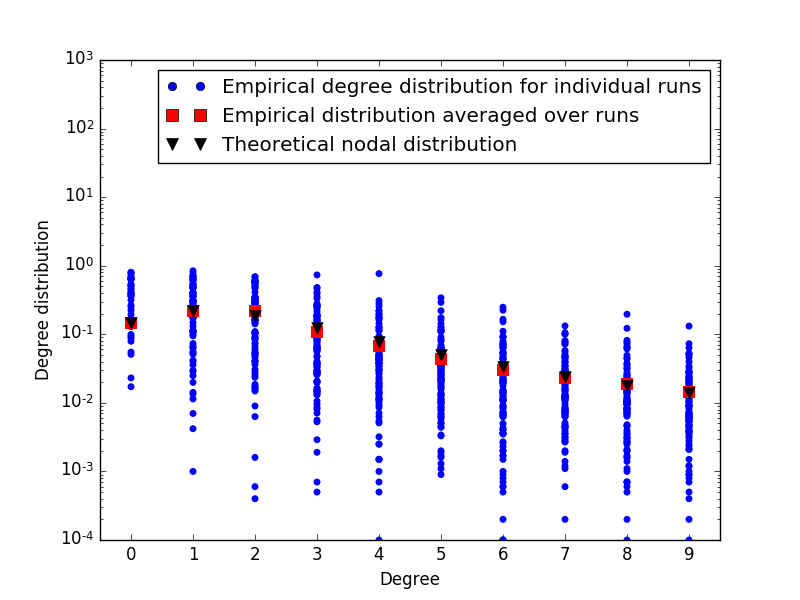}
 \caption{$n=10000$, $R=100$}
 \label{fig:Stripplot_10000nodes}
\end{subfigure}%
\begin{subfigure}{.55\textwidth}
 \centering
 \includegraphics[width=.8\linewidth]{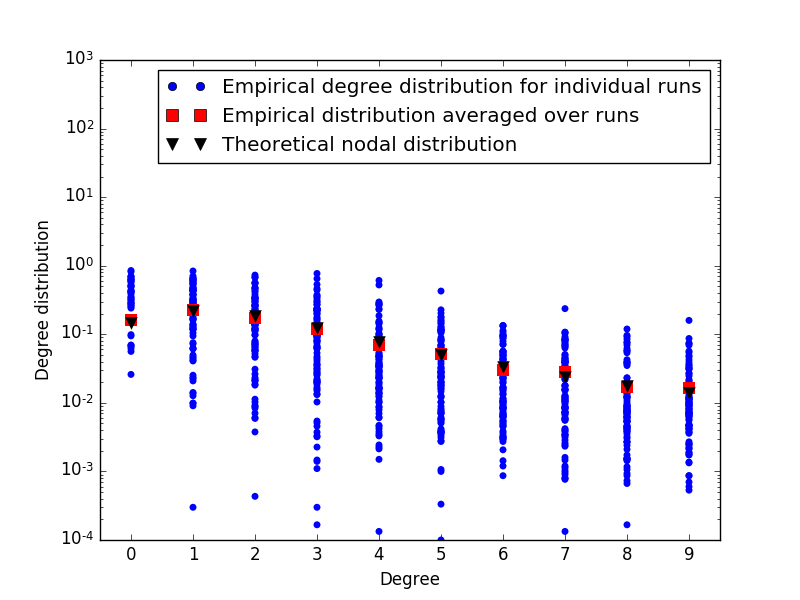}
 \caption{$n=30000$, $R=100$}
 \label{fig:Stripplot_30000nodes}
\end{subfigure}
\caption{The nodal degree distribution $p_{\rm Fuj}(.)$ was plotted against the  empirical degree distribution $\frac{N_n  ^{(r)}( . ;\theta _n ^{\star})}{n}$ for various runs $r=1,2,\ldots,R$. }
\label{fig:Stripplot_EmpiricalVsNodal}
\end{figure}

Next we explore the behavior of the empirical degree distribution (\ref{eq:EmpiricalPMF}) along the scaling
(\ref{eq:ExponentialRTgraphs_scaling}) (with $\lambda =1$) 
as generated through a single network realization.
We do so by plotting the histograms
\begin{equation}
\frac{N_n ^{(r)}(d;\theta_n ^\star)}{n} = \frac{1}{n} \sum_{k=1}^n \1{ D^{(r)}_{n,k} (\theta^\star_n) = d },
\quad 
\begin{array}{c}
d=0,1,\ldots \\
r=1, \ldots , R \\
\end{array}
\label{eq:SingleRealization}
\end{equation}
for various values of $d$ and $r$, and large $n$, and comparing against the corresponding value for $p_{\rm Fuj} (d)$.
In Figure~\ref{fig:Stripplot_EmpiricalVsNodal} we plot the histogram $\frac{N_n ^{(r)}( . ;\theta_n ^\star)}{n}$ for different runs $r=1,2,\ldots,R$ and varying graph sizes $n=10000,30000$, and observe high variability with respect to the nodal degree distribution $p_{\rm Fuj}(.)$, which does not change as the graph size is increased.

One might be tempted to smooth out the variability observed in Figure \ref{fig:Stripplot_EmpiricalVsNodal}
by averaging the empirical degree distributions (\ref{eq:SingleRealization}) over the $R$ i.i.d. realizations
$\mathbb{T}^{(1)}(n; \theta ^\star _n),  \mathbb{T}^{(2)}(n; \theta ^\star _n) ,  \ldots ,  \mathbb{T}^{(R)}(n; \theta ^\star _n)$,
resulting in the statistic
\begin{equation}
\frac{1}{R} \sum _{r=1} ^R \frac{N_n ^{(r)}(d;\theta_n ^\star)}{n}, 
\ d=0,1,\ldots .
\label{eq:Run_averaged_empirical_distribution}
\end{equation}
Fix $d=0,1, \ldots $. Under these circumstances, the Strong Law of Large Numbers yields
\begin{eqnarray}
\lim_{R \rightarrow \infty} \frac{1}{R} \sum _{r=1} ^R \frac{N_n ^{(r)}(d;\theta_n ^\star)}{n}
&=& \bE{ \frac{N_n (d;\theta_n ^\star)}{n} }
\quad \mbox{a.s.}
\end{eqnarray}
with
\[
 \bE{ \frac{N_n (d;\theta_n ^\star)}{n} }
 =
\bE{ \frac{1}{n} \sum_{k=1}^n \1{ D_{n,k} (\theta^\star_n) = d } }
= \bP{  D_{n} (\theta^\star_n) = d }
\]
by exchangeability. On the other hand, we have
$\lim_{n \rightarrow \infty} \bP{  D_{n} (\theta^\star_n) = d } = p_{\rm Fuj}(d)$ by virtue 
of Proposition \ref{prop:OneDimConvergence}.
Combining these observations yields the approximation
\begin{equation}
\frac{1}{R} \sum _{r=1} ^R \frac{N_n ^{(r)}(d;\theta_n ^\star)}{n} =_{\rm Approx} p_{\rm Fuj}(d)
\label{eq:Approx4}
\end{equation}
for large $n$ and $R$. The goodness of the approximation (\ref{eq:Approx4}) is noted in Figure~\ref{fig:Stripplot_EmpiricalVsNodal}, where the empirical distribution averaged over $R=100$ runs is observed to be very close to the nodal degree distribution.
However, the accuracy of the approximation (\ref{eq:Approx4}) does in no way imply the validity of  (\ref{eq:Intro_ConvergenceProbability}). 
In fact the mistaken belief that (\ref{eq:Intro_ConvergenceProbability}) holds,
implicitly assumed in the papers \cite{CCDM,SC}, might have stemmed from using the smoothed estimate (\ref{eq:Approx4}).


\section{A proof of Theorem \ref{thm:MainResult}}
\label{sec:ProofMainResult}

Fix $d=0,1, \ldots $. 
We establish the weak convergence of the sequence $\{ P_n(d;\theta^\star_n) , \ n=2,3, \ldots \}$
by arguments inspired by the method of moments; for details on the classical results concerning this approach,
see the references \cite[Thm. 4.5.5, p. 99]{ChungBook}  and \cite[Thm. 6.1, p. 140]{JansonLuczakRucinski}:

Proposition \ref{prop:ConvergenceMoments} suggests considering the 
mapping $\Phi_d: \mathbb{R} \rightarrow \mathbb{C}$ given by
\begin{equation}
\Phi _d (t)
= 1 + \sum_{r=1}^\infty \frac{(it)^r}{r!} \bP{ D_1 = d , \ldots , D_r = d },
\quad t \in \mathbb{R}.
\label{eq:Phi_d(t)_Defn}
\end{equation}
This definition is well posed with $\Phi_d(t)$ always an element of $\mathbb{C}$ since
\[
1 + \sum_{r=1}^\infty \frac{|t|^r}{r!} \bP{ D_1 = d , \ldots , D_r = d }
\leq 
 \sum_{r=0}^\infty \frac{|t|^r}{r!}  = e^{|t|},
\quad t  \in \mathbb{R}.
\]
In particular, the mapping $\Phi_d: \mathbb{R} \rightarrow \mathbb{C}$  is analytic on $\mathbb{R}$, hence continuous at $t=0$.
However, at this point in the proof, it is not yet known whether $\Phi_d$ is the characteristic function of a rv.

We close that gap as follows:
For each $n=2, 3, \ldots $, let
$\Phi_{d,n}: \mathbb{R} \rightarrow \mathbb{C}$ denote the characteristic function of the rv $P_n(d;\theta^\star_n)$,
i.e.,
\[
\Phi_{d,n} (t) = \bE{ e^{it P_n(d;\theta^\star_n) } },
\quad  t  \in \mathbb{R}.
\]
The obvious bound $0 \leq  P_n(d;\theta^\star_n)  \leq 1$ implies the uniform bounds
\[
\left |
1 + \sum_{r=1}^R \frac{(it)^r}{r!} P_n(d;\theta^\star_n)^r
\right |
\leq  1 + \sum_{r=1}^\infty \frac{|t|^r}{r!} P_n(d;\theta^\star_n)^r
\leq 
\sum_{r=0}^\infty \frac{|t|^r}{r!}  = e^{|t|},
\quad 
\begin{array}{c}
t  \in \mathbb{R}, \\
R=1,2, \ldots \\
\end{array}
\]
Therefore, applying the Bounded Convergence Theorem (with $R$ going to infinity), separately to the real and imaginary parts, 
we readily validate the series expansion 
\[
\Phi_{d,n} (t)  = 1 + \sum_{r=1}^\infty \frac{(it)^r}{r!} \bE{ P_n(d;\theta^\star_n)^r },
\quad  t  \in \mathbb{R}.
\]

For each $t $ in $\mathbb{R}$, it now follows that
\begin{eqnarray}
\Phi_{d,n} (t) - \Phi_{d} (t) 
&=&
 \sum_{r=1}^\infty \frac{(it)^r}{r!} 
 \left ( \bE{ \left ( P_n(d;\theta^\star_n) \right )^r } -  \bP{ D_1 = d , \ldots , D_r = d } \right ).
 \nonumber
\end{eqnarray}
Picking a positive integer $R$, we get

\begin{eqnarray}
\left | \Phi_{d,n} (t) - \Phi_{d} (t)  \right | 
\leq  \sum_{r=1}^R \frac{|t|^r}{r!} 
 \left |
 \bE{  P_n(d;\theta^\star_n)^r } -  \bP{ D_1 = d , \ldots , D_r = d } 
 \right |
 +  2 \sum_{r=R+1}^\infty \frac{|t|^r}{r!} .
 \nonumber 
\end{eqnarray}
For each $\varepsilon > 0$, there exists a positive integer $R^\star(\varepsilon,t)$ (independent of $n$) such that
\[
\sum_{r=R+1}^\infty \frac{|t|^r}{r!} 
\leq 
\varepsilon ,
\quad R \geq R^\star(\varepsilon,t),
\]
and on that range we obtain
\begin{eqnarray}
\lefteqn{ \limsup_{n \rightarrow \infty}  \left | \Phi_{d,n} (t) - \Phi_{d} (t)  \right | } & &
\nonumber \\
&\leq&
\sum_{r=1}^R  \frac{|t|^r}{r!} 
\limsup_{n \rightarrow \infty} 
\left (   \left |  \bE{  P_n(d;\theta^\star_n)^r } -  \bP{ D_1 = d , \ldots , D_r = d } \right | \right ) 
 +  2 \varepsilon
\nonumber
\end{eqnarray}
by the usual arguments.
Invoking Proposition \ref{prop:ConvergenceMoments} we readily conclude that
$\limsup_{n \rightarrow \infty}  \left | \Phi_{d,n} (t) - \Phi_{d} (t)  \right | \leq 2 \varepsilon$, whence
$ \lim_{n \rightarrow \infty} \Phi_{d,n} (t) = \Phi_{d} (t) $ since $\varepsilon > 0$ is arbitrary.

By a standard result due to Cram\'er and L\'evy 
\cite[Thm. 6.3.2, p. 161]{ChungBook}
\cite[Thm. 1, p. 320]{Shiryayev},
the mapping $\Phi_d: \mathbb{R} \rightarrow \mathbb{C}$, being continuous at $t=0$ with $\Phi_d(0)=1$,
must be the characteristic function of some rv, say $\Pi(d)$, namely
\begin{equation}
\bE{ e^{it \Pi(d)} } = \Phi_d(t),
\quad t \in \mathbb{R}
\label{eq:phi_d(t)}
\end{equation}
and the sequence $\{ P_n(d;\theta^\star_n) , \ n=2,3, \ldots \}$ converges weakly to the rv $\Pi(d)$.
\myendpf

For each $d=0,1, \ldots $, we read from (\ref{eq:Phi_d(t)_Defn}) and  (\ref{eq:phi_d(t)}) that $\bE{ \Pi(d) } = \bP{ D_1 = d}$ and
$\bE{ \Pi(d)^2 } = \bP{ D_1 = d, D_2=d}$, hence
\begin{eqnarray}
{\rm Var} \left [ \Pi(d) \right ] 
&=&  \bP{ D_1 = d , D_2 =d } - \bP{ D_1 = d } \bP{ D_2 = d } 
\nonumber \\
&=& {\rm Cov} \left [  \1{ D_1 = d } \1{ D_2 = d } \right ] .
\end{eqnarray}
The fact that ${\rm Var} \left [ \Pi(d) \right ]  > 0$ can be seen from the discussion at the end of Section \ref{sec:ProofMultDimConvergence}.
See also the references \cite{SPal_Thesis, PalMakowski-TNSE} for a discussion in the special case when the fitness rv is exponentially distributed.

\section{A proof of Proposition \ref{prop:ConvergenceMoments}}
\label{sec:ProofPropConvergenceMoments}

Fix $\theta > 0$, $d=0,1, \ldots $ and $r=1,2, \ldots $. 
For each $n=r,r+1, \ldots $, let
${\cal P}_{n,r}$ denote the collection of all {\em ordered} arrangements of $r$ {\em distinct} elements drawn from the set $\{ 1, \ldots , n \}$. 
Any such arrangement can be identified with a one-to-one mapping $\pi :\{1,  \ldots , r \} \rightarrow \{1,  \ldots , n \} $.

We begin with the well-known identity
\begin{eqnarray}
\prod_{s=0}^{r-1} \left ( N_n (d;\theta) - s \right )
=
\sum_{ \pi \in {\cal P}_{n,r} }  \left ( \prod_{t=1}^r \1{ D_{n,\pi(t)}(\theta) = d } \right ) .
\label{eq:Identity}
\end{eqnarray}
Taking expectations on both sides of (\ref{eq:Identity}), we obtain
\begin{eqnarray}
\bE{ \prod_{s=0}^{r-1} \left ( N_n (d;\theta) - s \right ) }
&=&
\sum_{\pi \in {\cal P}_{n,r}}  \bE{  \prod_{t=1}^r \1{ D_{n,\pi(t)}(\theta) = d }   }
\nonumber \\
&=&
\left | {\cal P}_{n,r} \right | 
\cdot \bP{ D_{n,1}(\theta) = d ,  \ldots , D_{n,r}(\theta) = d }  
\label{eq:ExpectationIdentity}
\end{eqnarray}
by the exchangeability of the rvs $D_{n,1}(\theta) , \ldots , D_{n,n}(\theta) $.
Dividing both sides of (\ref{eq:ExpectationIdentity}) by $n^r$, we get
\begin{eqnarray}
\bE{ \prod_{s=0}^{r-1} \left ( P_n (d;\theta) - \frac{s}{n} \right ) }
=
\frac{ \left | {\cal P}_{n,r} \right | }{n^r }
\cdot \bP{ D_{n,1}(\theta) = d ,  \ldots , D_{n,r}(\theta) = d }  
\label{eq:EQUALITY1}
\end{eqnarray}
with $ \left | {\cal P}_{n,r} \right | = n(n-1) \ldots (n-r+1)$.

Now consider the scaling $\theta^\star: \mathbb{N}_0 \rightarrow \mathbb{R}_+$ whose
existence is assumed in Assumption \ref{ass:A}.
For each $n=r, r+1, \ldots $ replace $\theta$ by $\theta^\star_n$ in (\ref{eq:EQUALITY1})
according to this scaling and let $n$ go to infinity in the resulting relation: 
Direct inspection shows that
$ \lim_{n \rightarrow \infty}  n^{-r} \left | {\cal P}_{n,r} \right | = 1$
while Proposition \ref{prop:MultiDimConvergence} yields
\[
\lim_{n \rightarrow \infty} 
\bP{ D_{n,1}(\theta^\star_n ) = d ,  \ldots , D_{n,r}(\theta^\star_n) = d }  
= \bP{ D_{1} = d ,  \ldots , D_{r} = d }
\]
where the rvs $D_1, \ldots , D_r$ are the limiting rvs appearing in the convergence (\ref{eq:MultiDimConvergence=r}).
Letting $n$ go to infinity in (\ref{eq:EQUALITY1}) yields
\begin{equation}
\lim_{n \rightarrow \infty} 
\bE{ \prod_{s=0}^{r-1} \left ( P_n (d;\theta^\star_n) - \frac{s}{n} \right ) }
=
\bP{ D_{1} = d ,  \ldots , D_{r} = d }.
\label{eq:EQUALITY2}
\end{equation}

Next, we readily check (say by recursion on $r$) that
\[
\left | \prod_{s=0}^{r-1} \left ( P_n (d;\theta^\star_n) - \frac{s}{n} \right ) - P_n (d;\theta^\star_n)^r \right |
\leq \sum_{s=0}^{r-1} \frac{s}{n}
= \frac{r(r-1)}{2n} ,
\quad n = r, r+1, \ldots 
\]
and it immediately follows that 
\[
\lim_{n \rightarrow \infty}
\bE{
\left | \prod_{s=0}^{r-1} \left ( P_n (d;\theta^\star_n) - \frac{s}{n} \right ) - P_n (d;\theta^\star_n)^r \right |
}
= 0.
\]
Combining this last fact with the convergence (\ref{eq:EQUALITY2}) yields
(\ref{eq:ConvergenceMoments=r}) by standard arguments, and the proof of Proposition \ref{prop:ConvergenceMoments} is now complete.
\myendpf

\section{A proof of Proposition \ref{prop:MultiDimConvergence}}
\label{sec:ProofMultDimConvergence}

The remainder of the paper is devoted to the proof of  Proposition \ref{prop:MultiDimConvergence}.
First some notation:
For each $r=1,2, \ldots$, 
let $\xi_{r|1} , \ldots , \xi_{r|r}$ denote the values of the
fitness rvs $\xi_1, \ldots , \xi_{r}$ arranged in increasing order, namely $\xi_{r|1} \leq  \ldots \leq \xi_{r|r}$,
with a lexicographic tiebreaker when needed. The rvs $\xi_{r|1} , \ldots , \xi_{r|r}$
are known as the {\em order statistics} associated with  the collection $\xi_1, \ldots , \xi_{r}$;
the rvs $\xi_{r|1} $ and $\xi_{r|r}$ are simply the minimum and maximum of the rvs 
$\xi_{1} , \ldots , \xi_{r}$, respectively \cite{DavidNagarajaBook}.
In what follows, the permutation $\alpha_r: \{1, \ldots , r \} \rightarrow \{1, \ldots, r \}$
arranges the rvs $\xi_1, \ldots , \xi_r$ in increasing order, i.e.,
\[
\xi_{r|s} = \xi_{\alpha_r(s)},
\quad s=1, \ldots , r
\]
(under the lexicographic tiebreaker) -- The permutation $\alpha_r$, being determined by the rvs
$\xi_{1} , \ldots , \xi_{r}$, it is a random permutation which is {\em uniformly} distributed over the group $\mathcal{S}_r$ 
of permutations of $\{1, \ldots , r \}$. Finally, with the notation introduced so far, write
\begin{equation}
G_r (z_1, \ldots , z_r)
=
\bE{
e^{ - \sum_{t=1}^r ( 1 - z_{\alpha_r (t)} ) \left ( \prod_{s=t+1}^r z_{\alpha_r(s)} \right ) \lambda(\xi_{r|t})  }
},
\quad 
\begin{array}{c}
0 \leq z_s \leq 1,\\
\quad s=1, \ldots , r.\\
\end{array}
\label{eq:G_r}
\end{equation}
By convention, the product of an empty set of factors is set to unity in the expression (\ref{eq:G_r}) and elsewhere in the discussion below.

The proof of Proposition \ref{prop:MultiDimConvergence} is an easy consequence of the following key analytical result.

\begin{proposition}
{\sl
Assume Assumption \ref{ass:A} to hold.
For each $r=1, 2, \ldots$,  we have
\begin{eqnarray}
\lim_{n \rightarrow \infty}
\bE{ \prod_{s=1}^r  z_s^{ D_{n,s}(\theta^\star_n) } }
= G_r (z_1, \ldots , z_r)
\label{eq:LimitPGFs}
\end{eqnarray}
for all $z_1, \ldots , z_r$ in $\mathbb{R}$ satisfying
\begin{equation}
0 \leq z_s \leq 1,
\quad s=1, \ldots , r.
\label{eq:RangeZZ}
\end{equation}
}
\label{prop:LimitPGFs}
\end{proposition}

This result is established in several steps which are presented from
Section \ref{sec:ProofPropLimitPGFsPart_I} to Section \ref{sec:ProofPropLimitPGFsPart_III}.
Proposition \ref{prop:LimitPGFs} does imply Proposition \ref{prop:MultiDimConvergence} by the usual
arguments: Indeed, by the Bounded Convergence
Theorem we get
\[
\lim_{ z_s \uparrow 1, \ s=1, \ldots , r} G_r (z_1, \ldots , z_r)
= 1
\]
and the mapping $G_r: [0,1]^r \rightarrow \mathbb{R}$ is therefore continuous at the point $(1, \ldots , 1)$.
This fact, coupled with the convergence (\ref{eq:LimitPGFs}), allows us to conclude that
$G_r$ is an $r$-dimensional pgf. Thus, there exists an $\mathbb{N}^r$-valued rv, denoted $(D_1, \ldots ,D_r)$, such that
\begin{equation}
\bE{ \prod_{s=1}^r  z_s^{ D_s } } = G_r(z_1, \ldots , z_r ),
\quad 
\begin{array}{c}
0 \leq z_s \leq 1,\\
\quad s=1, \ldots , r\\
\end{array}
\label{eq:LimitPGF2}
\end{equation}
and the convergence (\ref{eq:MultiDimConvergence=r})
follows in the usual manner.

For each $n=2,3, \ldots $, the rvs $D_{n,1}(\theta^\star_n) , \ldots , D_{n,r}(\theta^\star_n) $ are obviously exchangeable rvs,
and the exchangeability of the limiting rvs $D_1, \ldots ,D_r$ follows
because exchangeability is preserved under the weak convergence (\ref{eq:MultiDimConvergence=r}). 
This fact could also be gleaned directly from (\ref{eq:LimitPGF2}) as we note from (\ref{eq:G_r}) that
the mapping $G_r: [0,1]^r \rightarrow \mathbb{R}$ is permutation invariant in the sense that
\[
G_r(z_{\sigma(1)}, \ldots , z_{\sigma(r)})
=
G_r(z_1, \ldots , z_r ),
\quad 
\begin{array}{c}
0 \leq z_s \leq 1,\\
\quad s=1, \ldots , r\\
\end{array}
\]
for every permutation $\sigma$ of the index set $\{ 1, \ldots , r \}$: Indeed,
the random permutation $\sigma \circ \alpha_r: \{1, \ldots , r \} \rightarrow \{1, \ldots , r \}: s \rightarrow \sigma (\alpha_r (s))$
is uniform over $\mathcal{S}_r$ since the random permutation $\alpha_r: \{1, \ldots , r \} \rightarrow \{1, \ldots , r \}$ is itself uniform over $\mathcal{S}_r$,
and the rvs $\xi_1, \ldots , \xi_r$ are i.i.d. rvs; details are left to the interest reader.

Additional information can be obtained by direct inspection of  (\ref{eq:G_r}) and (\ref{eq:LimitPGF2}):
As expected, we retrieve  Proposition \ref{prop:OneDimConvergence} by looking at the case $r=1$, namely
\begin{equation}
\bE{ z_i^{ D_i}  } = \bE{ e^{-(1-z_i) \lambda(\xi_i) }},
\quad 
\begin{array}{c}
0 \leq z_i \leq 1\\
i= 1,2. \\
\end{array}
\label{eq:LimitPGF3a}
\end{equation}
For $r=2$, we also find 
\begin{equation}
\bE{ z_1^{ D_1} z_2^{ D_2}  } =
\bE{
e^{ -  ( 1 - z_{\alpha_2 (1)} )z_{\alpha_2(2)} \lambda(\xi_{2|1})  - (1-z_{\alpha_2(2)}) \lambda(\xi_{2|2}) }
},
\quad 
0 \leq z_1,z_2 \leq 1.
\label{eq:LimitPGF3b}
\end{equation}
Comparing (\ref{eq:LimitPGF3a}) and (\ref{eq:LimitPGF3b}) we can check that 
\begin{equation}
\bE{ z_1^{ D_1} z_2^{ D_2}  } 
\neq 
\bE{ z_1^{ D_1} } \bE{ z_2^{ D_2}  },
\quad 
0 \leq z_1,z_2 \leq 1
\end{equation}
and the rvs $D_1, \ldots , D_r$ are therefore not independent.
This completes the proof of  Proposition \ref{prop:MultiDimConvergence}.
\myendpf

To further illustrate this last point, consider the special case when the mapping $\lambda : \mathbb{R}_+ \rightarrow \mathbb{R}_+$
appearing in Assumption \ref{ass:A} is constant, say $\lambda (x) = c$ for all $x \geq 0$ with $c>0$, as would be the case 
for the Pareto distribution (\ref{eq:ParetoFitness}).
The expressions (\ref{eq:LimitPGF3a}) and (\ref{eq:LimitPGF3b}) now become
\begin{equation}
\bE{ z_i^{ D_i}  } = e^{-(1-z_i) c }
\quad 
\begin{array}{c}
0 \leq z_i \leq 1\\
i= 1,2 \\
\end{array}
\label{eq:LimitPGF3A}
\end{equation}
and
\begin{equation}
\bE{ z_1^{ D_1} z_2^{ D_2}  } = e^{ -  ( 1 - z_1z_2) c },
\quad 
0 \leq z_1,z_2 \leq 1.
\label{eq:LimitPGF3B}
\end{equation}
Thus, each of the rvs $D_1, \ldots , D_r$ is Poisson distributed with parameter $c$ and $\bP{ D_1 = D_2 } = 1$ since
\[
\bP{ D_1 =d, D_2 = d^\prime } 
= \frac{ c^d}{d!} e^{-c } \cdot  \delta (d,d^\prime ),
\quad d, d^\prime =01,2, \ldots 
\]
Exchangeability yields $\bP{ D_1 = D_j } = 1$ for every $j=2, \ldots , r$, hence
$\bP{ D_1 = \ldots = D_r } = 1$, and the 
rvs $D_1, \ldots , D_r$  are certainly not independent! Moreover, for each $d=0,1, \ldots $, we get 
$\bP{ D_1 = d, \ldots , D_r = d } = \bP{ D_1 =d }$ for all $r=1,2, \ldots $, whence
\begin{eqnarray}
\Phi_d(t)
&=& 
1 + \sum_{r=1}^\infty \frac{ (it)^r }{r!} \bP{ D_1 =d } 
\nonumber \\
&=&
1 + \left ( e^{it} - 1 \right ) \cdot \bP{ D_1 =d }
\nonumber \\
&=&
\left ( 1 - \bP{ D_1 = d } \right ) + e^{it} \cdot \bP{ D_1 =d },
\quad t \in \mathbb{R}.
\end{eqnarray}
In other words, the distribution of the rv $\Pi(d)$ is the two point mass distribution $( 1 - \bP{ D_1 = d } ,  \bP{ D_1 =d } )$
on the set $\{ 0, 1 \}$ with $\bP{ D_1 = d } = \frac{c^d}{d!} e^{-c}$.
Obviously, 
\[
{\rm Var} [ \Pi (d) ] = \bP{ D_1 = d } - \bP{ D_1 = d } \bP{ D_1 = d } =  \bP{ D_1 = d } \left ( 1 - \bP{ D_1 = d } \right )   > 0 
\]
since $0 <  \bP{ D_1 = d } < 1 $.

When the fitness rv $\xi$ is exponentially distributed, 
explicit expressions were obtained for ${\rm Var} [ \Pi (d) ] $ by direct arguments in the earlier references
\cite{SPal_Thesis, PalMakowski_CDC2015, PalMakowski-TNSE}.

\section{A proof of Proposition \ref{prop:LimitPGFs} -- A reduction step}
\label{sec:ProofPropLimitPGFsPart_I} 

Throughout this section the integer $r=1,2, \ldots $  and the parameter $\theta > 0$ are held fixed.
Pick $n > r$. For each $k=1, \ldots , r$, we write
\[
D_{n,k} (\theta) 
= \sum_{\ell=1, \ \ell \neq k}^r \1{ \xi_k + \xi_\ell  > \theta }
+
D^{(r)}_{n,k} (\theta)
\]
with
\[
D^{(r)}_{n,k} (\theta)
= \sum_{\ell=r+1}^n \1{ \xi_k + \xi_\ell  > \theta }.
\]

As the scaling $\theta^\star : \mathbb{N}_0 \rightarrow \mathbb{R}_+$ satisfies
$\lim_{n \rightarrow \infty} \theta^\star_n = \infty$, it is plain that
\[
\lim_{n \rightarrow \infty} 
\max
\left ( 
\sum_{\ell=1, \ \ell \neq k}^r \1{ \xi_k + \xi_\ell  > \theta^\star_n },
\ k=1, \ldots , r 
\right ) = 0
\quad a.s.
\]
and (\ref{eq:LimitPGFs}) takes place if and only if
\begin{equation}
\lim_{n \rightarrow \infty} 
\bE{ \prod_{s=1}^r z_s^{ D^{(r)}_{n,s} (\theta^\star_n) } }
=
G_r (z_1, \ldots , z_r)
\label{eq:LimitJointPGFs_0}
\end{equation}
for all $z_1, \ldots , z_r$ in $\mathbb{R}$ which satisfy (\ref{eq:RangeZZ}).

Our first step towards establishing (\ref{eq:LimitJointPGFs_0}) is to evaluate the joint pgfs.
Pick $z_1, \ldots , z_r$ in $\mathbb{R}$.
Under the enforced independence assumptions, it is plain that
\begin{eqnarray}
\bE{ \prod_{s=1}^r z_s^{ D^{(r)}_{n,s} (\theta) } }
&=&
\bE{ \prod_{s=1}^r z_s^{  \sum_{\ell=r+1}^n \1{ \xi_s + \xi_\ell  > \theta } } }
\nonumber \\
&=&
\bE{ \prod_{s=1}^r \prod_{\ell=r+1}^n z_s^{ \1{ \xi_s + \xi_\ell  > \theta } } }
\nonumber \\
&=&
\bE{ \prod_{\ell=r+1}^n \prod_{s=1}^r  z_s^{ \1{ \xi_s + \xi_\ell  > \theta } } }
\nonumber \\
&=&
\bE{ 
\bE{ \prod_{\ell=r+1}^n \prod_{s=1}^r  z_s^{ \1{ \xi_s + \xi_\ell  > \theta } } \Bigl | \xi_1, \ldots , \xi_r }
}
\nonumber \\
&=&
\bE{ 
\bE{ \prod_{\ell=r+1}^n \prod_{s=1}^r  z_s^{ \1{ x_s + \xi_\ell  > \theta } } }_{ x_1 = \xi_1, \ldots , x_r=\xi_r}
}.
\label{eq:Reduction2}
\end{eqnarray}

With arbitrary $x_1, \ldots , x_r$ in $\mathbb{R}_+$, we get
\begin{eqnarray}
\bE{ \prod_{\ell=r+1}^n \prod_{s=1}^r  z_s^{ \1{ x_s + \xi_\ell  > \theta } } }
&=&
\prod_{\ell=r+1}^n \bE{ \prod_{s=1}^r  z_s^{ \1{ x_s + \xi_\ell  > \theta } } }
\nonumber \\
&=&
F_r (\theta; z_1, \ldots , z_r; x_1, \ldots , x_r)^{n-r}
\nonumber
\end{eqnarray}
where we have set
\begin{eqnarray}
F_r (\theta; z_1, \ldots , z_r; x_1, \ldots , x_r)
&=& \bE{ \prod_{s=1}^r  z_s^{ \1{ x_s + \xi > \theta } } } 
\nonumber \\
&=&
\bE{ \prod_{s=1}^r  \left (  \1{ x_s + \xi > \theta }  z_s + \1{ x_s + \xi \leq \theta }  \right ) }
\nonumber \\
&=&
\bE{ \prod_{s=1}^r  \left (  1 - \left ( 1 - z_s \right ) \1{ x_s + \xi > \theta } \right ) }.
\label{eq:F_r}
\end{eqnarray}
Substituting back into (\ref{eq:Reduction2}) we conclude that
\begin{equation}
\bE{ \prod_{s=1}^r z_s^{ D^{(r)}_{n,s} (\theta) } }
=
\bE{ F_r (\theta; z_1, \ldots , z_r; \xi_1 , \ldots , \xi_r )^{n-r} },
\label{eq:JointPGF}
\end{equation}
and the desired result (\ref{eq:LimitJointPGFs_0}) does hold if we show that
\begin{equation}
\lim_{n \rightarrow \infty} 
\bE{ F_r (\theta^\star_n; z_1, \ldots , z_r; \xi_1 , \ldots , \xi_r )^{n-r} }
=
G_r (z_1, \ldots , z_r)
\label{eq:LimitJointPGF}
\end{equation}
for all $z_1, \ldots , z_r$ in $\mathbb{R}$ which satisfy (\ref{eq:RangeZZ}).

\section{A proof of Proposition \ref{prop:LimitPGFs} -- A decomposition}
\label{sec:ProofPropLimitPGFsPart_II} 

Fix $\theta > 0$ and $r=1,2, \ldots $.
To further analyze the expression (\ref{eq:JointPGF}) with $x_1, \ldots , x_r$ in $\mathbb{R}_+$, we introduce the index set
\[
{\mathcal S}(\theta; x_1, \ldots , x_r)
=
\left \{
s=1, \ldots , r: \ x_s > \theta 
\right \}.
\]
There are two possibilities which we now explore in turn: 
Either ${\mathcal S}(\theta; x_1, \ldots , x_r)$ is empty or it is not, leading to a natural decomposition
expressed through Lemmas \ref{lem:Decomposition=A} and \ref{lem:Decomposition=B}.

\begin{lemma}
{\sl With $x_1, \ldots , x_r$ in $\mathbb{R}_+$, whenever 
${\mathcal S}(\theta; x_1, \ldots , x_r)$ is non-empty, we have
\begin{eqnarray}
\lefteqn{
F_r (\theta; z_1, \ldots , z_r; x_1 , \ldots , x_r )
} & &
\nonumber \\
&=&
\left ( \prod_{s \in {\mathcal S}(\theta; x_1, \ldots , x_r) } z_s  \right )
\cdot
\bE{ \prod_{s \notin {\mathcal S}(\theta; x_1, \ldots , x_r)}^r  \left (   1 - \left ( 1-z_s \right ) \1{ x_s + \xi > \theta }  \right ) }
\label{eq:Decomposition=A}
\end{eqnarray}
for all $z_1, \ldots , z_r$ in $\mathbb{R}$.
}
\label{lem:Decomposition=A}
\end{lemma}

\myproof
Pick arbitrary $x_1, \ldots , x_r$ in $\mathbb{R}_+$ with non-empty ${\mathcal S}(\theta; x_1, \ldots , x_r)$.
For all $z_1, \ldots , z_r$ in $\mathbb{R}$, it is easy to check by direct inspection from the expression
(\ref{eq:F_r}) that (\ref{eq:Decomposition=A}) holds since 
$1 - \left ( 1-z_s \right ) \1{ x_s + \xi > \theta }   = z_s$
whenever $s$ belongs to ${\mathcal S}(\theta; x_1, \ldots , x_r)$.
\myendpf

As an immediate consequence of (\ref{eq:Decomposition=A})  we have the inequality
\begin{eqnarray}
0 \leq F_r (\theta; z_1, \ldots , z_r; x_1 , \ldots , x_r ) \leq 1,
\quad
\begin{array}{c}
x_1, \ldots , x_r \in \mathbb{R}_+ \\
\mbox{with} \\
|{\mathcal S}(\theta; x_1, \ldots , x_r)| > 0 \\
\end{array}
\label{eq:DecompositionInequality=1}
\end{eqnarray}
for all $z_1, \ldots , z_r$ in $\mathbb{R}$ in the range
\begin{equation}
|z_s| \leq 1, \quad s=1, \ldots , r.
\label{eq:RangeZ}
\end{equation}
This is because,  it is always the case there that
$\left | 1 - \left ( 1-z_s \right ) \1{ x_s + \xi > \theta }   \right | \leq 1 $ if $|z_s| \leq 1$.

We now turn to the case when the index set ${ \mathcal S}(\theta; x_1, \ldots , x_r)$ is empty,
a fact characterized by the conditions
\begin{equation}
x_s \leq \theta ,
\quad s=1, \ldots , r.
\label{eq:IndexSetEmpty}
\end{equation}
It will be convenient to arrange the values $x_1, \ldots , x_r$ in increasing order, say 
$ x_{(r|1)} \leq x_{(r|2)} \leq \ldots \leq x_{(r|r)}$, with a lexicographic tiebreaker.
Let $a_r$ be any permutation of $\{ 1, \ldots , r \}$
such that $x_{(r|s)} = x_{a_r(s)} $ for all $s=1, \ldots , r$ -- Obviously this permutation is determined by the values  $x_1, \ldots , x_r$.
In what follows we shall use the convention $x_{(r|0)} = -\infty$ and $ x_{(r|r+1)} = \infty$. 

\begin{lemma}
{\sl 
With $x_1, \ldots , x_r$ in $\mathbb{R}_+$, whenever 
${\mathcal S}(\theta; x_1, \ldots , x_r)$ is empty, we have
\begin{eqnarray}
F_r(\theta ; z_1, \ldots , z_r ; x_1, \ldots , x_r )
=
\sum_{t=0}^{r} \left ( \prod_{s=t+1}^r z_{a_r (s)} \right ) \cdot 
\left ( F( \theta - x_{(r|t)} ) - F( \theta - x_{(r|t+1)} )
\right )
\label{eq:Decomposition=B}
\end{eqnarray}
for all $z_1, \ldots , z_r$ in $\mathbb{R}$.
}
\label{lem:Decomposition=B}
\end{lemma}

\myproof
In what follows, the values $z_1, \ldots , z_r$ in $\mathbb{R}$ are held fixed.
Given $x_1, \ldots , x_r$ in $\mathbb{R}_+$ and $\theta > 0$, we define the events
\[
A_{r|t} ( x_1, \ldots , x_r ; \theta) = 
\left [ x_{(r|t)} + \xi \leq \theta < x_{(r|t+1)} + \xi \right ],
\quad t=0,1,, \ldots , r.
\]
Under the enforced conventions, we have
$A_{r|0} ( x_1, \ldots , x_r ; \theta) = [ \theta < x_{(r|1)} + \xi ]$ and $A_{r|r} ( x_1, \ldots , x_r ; \theta) = [ x_{(r|r)} + \xi \leq \theta ]$.
When ${\mathcal S}(\theta; x_1, \ldots , x_r)$ is empty, the $r+1$ events $A_{r|0} ( x_1, \ldots , x_r ; \theta), \ldots , A_{r|r} ( x_1, \ldots , x_r ; \theta)$ 
are pairwise disjoint and form a partition of the sample space. Using this fact in the expression (\ref{eq:F_r}) we find
\begin{eqnarray}
\lefteqn{
F_r(\theta; z_1, \ldots , z_r; x_1, \ldots , x_r)
} & &
\nonumber \\
&=&
\sum_{t=0}^r 
\bE{ 
\1{ A_{r|t} (x_1, \ldots , x_r ; \theta)  }  \prod_{s=1}^r  \left (  1 - \left ( 1 - z_s \right ) \1{ x_s + \xi > \theta } \right ) 
}.
\label{eq:Decomposition}
\end{eqnarray}

(i) On the event $A_{r|0} ( x_1, \ldots , x_r ; \theta)$, we have
$\theta < x_{(r|1)} + \xi $, thus $ \theta < x_{s} + \xi $ for all $s =1, \ldots , r$, so that
\[
\prod_{s=1}^r  \left (   1 - \left ( 1-z_s \right ) \1{ x_s + \xi > \theta }  \right ) 
=
\prod_{s=1}^r z_s ,
\]
whence
\begin{eqnarray}
\lefteqn{ \bE{ \1{  A_{r|0} ( x_1, \ldots , x_r ; \theta) }
\prod_{s=1}^r  \left (   1 - \left ( 1-z_s \right ) \1{ x_s + \xi > \theta }  \right ) }
} & &
\nonumber \\
&=&
\left ( \prod_{s=1}^r z_s  \right ) \cdot \bP{ \theta < x_{(r|1)} + \xi }
\nonumber \\
&=&
\left ( \prod_{s=1}^r z_s  \right ) \cdot  \left ( 1 - F( \theta - x_{(r|1)} ) \right ) .
\label{eq:DecompositionFirstTerm}
\end{eqnarray}

(ii) With $t=1, \ldots , r-1$,  on the event $A_{r|t} ( x_1, \ldots , x_r ; \theta) $ it holds that
$ x_{(r|1)} + \xi \leq \theta, \ldots , x_{(r|t)} + \xi  \leq \theta $ and $\theta < x_{(r|t+1)} + \xi , \ldots , \theta <  x_{(r|r)} + \xi $, whence
\[
\prod_{s=1}^r  \left (   1 - \left ( 1-z_s \right ) \1{ x_s + \xi > \theta }  \right )
= 
\left ( \prod_{s=t+1}^r z_{a_r(s)} \right ) .
\]
We readily conclude to
\begin{eqnarray}
\lefteqn{
\bE{ \1{  A_{r|t} ( x_1, \ldots , x_r ; \theta) }
\prod_{s=1}^r  \left (   1 - \left ( 1-z_s \right ) \1{ x_s + \xi > \theta }  \right ) }
}
& &
\nonumber \\
&=&
\left ( \prod_{s=t+1}^r z_{a_r(s)} \right )
\cdot
\bP{  x_{(r|t)} + \xi \leq \theta < x_{(r|t+1)} + \xi }
\nonumber \\
&=&
\left ( \prod_{s=t+1}^r z_{a_r(s)} \right )
\cdot
\left (
F( \theta - x_{(r|t)} ) - F( \theta - x_{(r|t+1)} )
\right ).
\label{eq:DecompositionSecondTerm}
\end{eqnarray}

(iii) Finally, on the event $A_{r|r} ( x_1, \ldots , x_r ; \theta)$, it holds that
$ x_{(r|r)} + \xi \leq \theta $, thus $ x_{s} + \xi \leq \theta$ for all $s =1, \ldots , r$, so that
\[
\prod_{s=1}^r  \left (   1 - \left ( 1-z_s \right ) \1{ x_s + \xi > \theta }  \right )  = 1,
\]
whence
\begin{eqnarray} 
\bE{ \1{  A_{r|r} ( x_1, \ldots , x_r ; \theta) }
\prod_{s=1}^r  \left (   1 - \left ( 1-z_s \right ) \1{ x_s + \xi > \theta }  \right ) }
&=&
\bP{ x_{(r|r)} + \xi \leq \theta }
\nonumber \\
&=&
F( \theta - x_{(r|r)} )  .
\label{eq:DecompositionThirdTerm}
\end{eqnarray}

To complete the proof we substitute 
(\ref{eq:DecompositionFirstTerm}), (\ref{eq:DecompositionSecondTerm}) and (\ref{eq:DecompositionThirdTerm})
into (\ref{eq:Decomposition}), and recall  that $F( \theta - x_{(r|0)} ) = 1$ and $F( \theta - x_{(r|r+1)} ) = 0$ 
under the conventions adopted here.
\myendpf

\section{A proof of Proposition \ref{prop:LimitPGFs} -- Taking the limit}
\label{sec:ProofPropLimitPGFsPart_III} 

In order to establish the convergence (\ref{eq:LimitJointPGFs_0}) we return to the
expression (\ref{eq:JointPGF}) for the joint pgf of the relevant rvs.

\subsection{ A useful intermediary fact}

Fix $n =2,3, \ldots$ with $r < n$.
For arbitrary $\theta > 0$, consider $x_1, \ldots , x_r$ in $\mathbb{R}_+$ and $z_1, \ldots , z_r$ in $\mathbb{R}$.
In what follows it will be convenient to define
\[
\Lambda_r (\theta; z_1, \ldots , z_r;x_1, \ldots , x_r)
= 1 - F_r (\theta; z_1, \ldots , z_r;x_1, \ldots , x_r)
\]
so that
\[
F_r (\theta; z_1, \ldots , z_r;x_1, \ldots , x_r)
=
1 - \Lambda_r (\theta; z_1, \ldots , z_r;x_1, \ldots , x_r).
\]

Whenever ${\cal S}(\theta; x_1, \ldots , x_r)$ is empty, Lemma \ref{lem:Decomposition=B} gives
\begin{eqnarray}
\lefteqn{
\Lambda_r (\theta; z_1, \ldots , z_r;x_1, \ldots , x_r)
} & & 
\nonumber \\
&=&
1 - \sum_{t=0}^{r} \left ( \prod_{s=t+1}^r z_{a_r (s)} \right ) \cdot 
\left ( F( \theta - x_{(r|t)} ) - F( \theta - x_{(r|t+1)} )
\right )
\nonumber \\
&=&
- \sum_{t=0}^{r-1} \left ( \prod_{s=t+1}^r z_{a_r (s)} \right ) \cdot 
\left ( F( \theta - x_{(r|t)} ) - F( \theta - x_{(r|t+1)})  \right )
+ \left ( 1 -   F( \theta - x_{(r|r)} ) \right ) .
\label{eq:Lambda}
\end{eqnarray}

Replace $\theta$ by $\theta^\star_n$ in (\ref{eq:Lambda})
according to the scaling $\theta^\star : \mathbb{N}_0 \rightarrow \mathbb{R}_+$ 
stipulated in Assumption \ref{ass:A}: Letting $n$ go to infinity in the resulting relation, we get
\[
\lim_{n \rightarrow \infty} 
n \left ( 1 -   F( \theta^\star_n - x_{(r|r)} ) \right ) = \lambda ( x_{(r|r)} ) 
\]
and
\begin{eqnarray}
\lefteqn{
\lim_{n \rightarrow \infty} 
n \left ( F( \theta^\star_n - x_{(r|t)} ) - F( \theta^\star_n - x_{(r|t+1)} ) \right )
} & &
\nonumber \\
&=&
 \lim_{n \rightarrow \infty} 
n \left ( 
\left ( 1 - F( \theta^\star_n - x_{(r|t+1)} \right ) - \left ( 1 - F( \theta^\star_n - x_{(r|t)} ) \right )
\right )
\nonumber \\
&=&
\left \{
\begin{array}{ll}
\lambda ( x_{(r|1)} ) & \mbox{if $t=0$} \\
  & \\
\lambda ( x_{(r|t+1)} )  -  \lambda ( x_{(r|t)} ) & \mbox{if $t=1, \ldots , r-1$.} \\
\end{array}
\right .
\end{eqnarray}
As a result, with ${\cal S}(\theta; x_1, \ldots , x_r)$ empty, we have
\begin{eqnarray}
\lefteqn{
\lim_{n \rightarrow \infty}  n \Lambda_r (\theta^\star_n; z_1, \ldots , z_r;x_1, \ldots , x_r)
} & &
\nonumber \\
&=&
- \lambda ( x_{(r|1)} ) \left ( \prod_{s=1}^r z_s \right )  - \sum_{t=1}^{r-1} 
\left ( \prod_{s=t+1}^r z_{a_r (r|s)} \right ) \left ( \lambda ( x_{(r|t+1)} )  - \lambda ( x_{(r|t)} ) \right ) +  \lambda ( x_{(r|r)} )
\nonumber \\
&=&
- \sum_{t=1}^r \lambda(x_{(r|t)})
\left ( \prod_{s=t}^r z_{a_r(r|s)} - \prod_{s=t+1}^r z_{a_r(r|s)} \right )
\nonumber \\
&=&
\sum_{t=1}^r \lambda(x_{(r|t)}) ( 1 - z_{a_r (t)} ) \prod_{s=t+1}^r z_{a_r(s)} ,
\end{eqnarray}
and the conclusion
\begin{eqnarray}
\lim_{n \rightarrow \infty} F_r (\theta^\star_n; z_1, \ldots , z_r;x_1, \ldots , x_r)^{n-r}
&=&
\lim_{n \rightarrow \infty} 
\left ( 1 - \Lambda_r (\theta^\star_n; z_1, \ldots , z_r;x_1, \ldots , x_r )  \right )^{n-r}
\nonumber \\
&=&
e^{ - \sum_{t=1}^r \lambda(x_{(r|t)}) ( 1 - z_{a_r (t)} ) \prod_{s=t+1}^r z_{a_r(s)}  }
\nonumber
\end{eqnarray}
follows by standard arguments \cite[Prop. 3.1.1., p. 116]{EKM}.

\subsection{In the limit}

Pick $z_1, \ldots , z_r$ in $\mathbb{R}$ such that (\ref{eq:RangeZZ}) holds.
For each $ n=2,3, \ldots $ with $r < n$, the decomposition
\begin{eqnarray}
\lefteqn{
\bE{ F_r (\theta^\star_n; z_1, \ldots , z_r; \xi_1 , \ldots , \xi_r )^{n-r} } 
} & &
\nonumber \\
&=&
\bE{ \1{ | {\mathcal S}(\theta^\star_n; \xi_1, \ldots , \xi_r) |  > 0 }  F_r (\theta^\star_n; z_1, \ldots , z_r; \xi_1 , \ldots , \xi_r )^{n-r} } 
\nonumber \\
& &
~+
\bE{ \1{ | {\mathcal S}(\theta^\star_n; \xi_1, \ldots , \xi_r) |  = 0 }  F_r (\theta^\star_n; z_1, \ldots , z_r; \xi_1 , \ldots , \xi_r )^{n-r} } 
\label{eq:DECOMPOSITION}
\end{eqnarray}
holds.

Because $\lim_{n \rightarrow \infty} \theta^\star_n = \infty$, it follows that
\[
\lim_{n \rightarrow \infty}
\1{   | {\mathcal S}(\theta^\star_n; \xi_1, \ldots , \xi_r) |  = 0  }
=
\lim_{n \rightarrow \infty}
\1{ \xi_1 \leq \theta^\star_n , \ldots ,  \xi_r \leq \theta^\star_n }
= 1,
\]
whence
$\lim_{n \rightarrow \infty} \bP{   | {\mathcal S}(\theta^\star_n; \xi_1, \ldots , \xi_r) |  = 0  } = 1 $.
The inequality (\ref{eq:DecompositionInequality=1}) now implies
\begin{equation}
\lim_{n \rightarrow \infty}
\bE{ \1{ | {\mathcal S}(\theta^\star_n; \xi_1, \ldots , \xi_r) |  > 0 }  F_r (\theta^\star_n; z_1, \ldots , z_r; \xi_1 , \ldots , \xi_r )^{n-r} }  
= 0
\label{eq:LimitFirstTerm}
\end{equation}
since the condition  (\ref{eq:RangeZZ}) is more restrictive than (\ref{eq:RangeZ}).

Next, for each $ n=2,3, \ldots $ with $r < n$, we see that
\begin{eqnarray}
\lefteqn{
\bE{ \1{ | {\mathcal S}(\theta^\star_n; \xi_1, \ldots , \xi_r) |  = 0 }  F_r (\theta^\star_n; z_1, \ldots , z_r; \xi_1 , \ldots , \xi_r )^{n-r} } 
} & &
\nonumber \\
&=&
\bE{ \1{ | {\mathcal S}(\theta^\star_n; \xi_1, \ldots , \xi_r) |  = 0 }  
\left ( 1 - \Lambda_r (\theta^\star_n; z_1, \ldots , z_r; \xi_1 , \ldots , \xi_r ) \right ) ^{n-r} } 
\end{eqnarray}
with (\ref{eq:Lambda}) yielding
\begin{eqnarray}
\lefteqn{
\Lambda_r (\theta^\star_n; z_1, \ldots , z_r;\xi_1, \ldots ,\xi_r)
} & & 
\nonumber \\
&=&
1 -   F( \theta^\star_n - \xi_{r|r} )
- \sum_{t=0}^{r-1} 
\left ( \prod_{s=t+1}^r z_{\alpha_r (s)} \right )
\cdot 
\left ( F( \theta^\star_n - \xi_{r|t} ) - F( \theta^\star_n  - \xi_{r|t+1}  ) \right ).
\label{eq:Lambda2}
\end{eqnarray}
As we pass from (\ref{eq:Lambda}) to (\ref{eq:Lambda2}), we recall that the order statistics $ \xi_{r|1} , \ldots , \xi_{r|r}$ 
associated with $\xi_1, \ldots , \xi_r$ were introduced
in the statement of Proposition \ref{prop:LimitPGFs}, together with the {\em random} permutation 
$\alpha_r: \{ 1, \ldots, r \} \rightarrow \{ 1, \ldots, r \} $.  The random permutation $\alpha_r$ coincides with the deterministic permutation
$a_r$ induced by the values $x_1, \ldots , x_r$ with $x_1 = \xi_1, \ldots , x_r = \xi_r$.

Under the condition (\ref{eq:RangeZZ}) it is plain that
\[
0 \leq \left ( \prod_{s=t+1}^r z_{\alpha_r (s)} \right )   \leq 1,
\quad t=0, \ldots , r-1,
\]
while the bounds
$0 \leq \Lambda_r (\theta^\star_n; z_1, \ldots , z_r;\xi_1, \ldots ,\xi_r) \leq 1 $ hold by direct inspection of (\ref{eq:Lambda2}), whence
$\left |
F_r (\theta^\star_n; z_1, \ldots , z_r;\xi_1, \ldots ,\xi_r)^{n-r}
\right |  \leq 1$.
With the fact
$ \lim_{n \rightarrow \infty} \1{ | {\mathcal S}(\theta^\star_n; \xi_1, \ldots , \xi_r) |  = 0 }  = 1$
noted earlier, we see that the  convergence
\begin{eqnarray}
\lefteqn{
\lim_{n \rightarrow \infty}
\1{ | {\mathcal S}(\theta^\star_n; \xi_1, \ldots , \xi_r) |  = 0 }  F_r (\theta^\star_n; z_1, \ldots , z_r; \xi_1 , \ldots , \xi_r )^{n-r} 
} & &
\nonumber \\
&=&
e^{ - \sum_{t=1}^r \lambda(\xi_{r|t}) ( 1 - z_{\alpha_r (t)} ) \prod_{s=t+1}^r z_{\alpha_r(s)}  }
\nonumber
\end{eqnarray}
takes place boundedly, and the Bounded Convergence  Theorem can then be applied to yield
\begin{eqnarray}
\lefteqn{
\lim_{n \rightarrow \infty}
\bE{
\1{ | {\mathcal S}(\theta^\star_n; \xi_1, \ldots , \xi_r) |  = 0 }  F_r (\theta^\star_n; z_1, \ldots , z_r; \xi_1 , \ldots , \xi_r )^{n-r} 
}
} & &
\nonumber \\
&=&
\bE{ e^{- \sum_{t=1}^r \lambda(\xi_{r|t}) ( 1 - z_{\alpha (t)} ) \prod_{s=t+1}^r z_{\alpha_r(s)}  } }
\nonumber \\
&=& G_r (z_1, \ldots , z_r) .
\label{eq:LimitSecondTerm}
\end{eqnarray}

Let $n$ go to infinity in (\ref{eq:DECOMPOSITION}):
Collecting (\ref{eq:LimitFirstTerm}) and (\ref{eq:LimitSecondTerm}), and using 
(\ref{eq:DECOMPOSITION}) we conclude that
(\ref{eq:LimitJointPGF}) indeed holds on the range (\ref{eq:RangeZZ}).
\myendpf

\section*{Acknowledgment}
This work was supported by NSF Grant CCF-1217997.
The paper was completed during the academic year 2014-2015 while A.M. Makowski 
was a Visiting Professor with the Department of Statistics of the Hebrew University of Jerusalem 
with the support of a fellowship from the Lady Davis Trust.


\begin{thebibliography}{1}

\bibitem{BarabasiAlbert}
A.-L. Barab\'{a}si and R. Albert,
\lq\lq Emergence of scaling in random networks,"
{\sl Science} {\bf 286} (1999), pp. 509-512.


\bibitem{BillingsleyBook}
P. Billingsley,
{\sl Convergence of Probability Measures},
John Wiley \& Sons, New York (NY), 1968.



\bibitem{Bollobas_BAmodel}
B. Bollob\'{a}s, O. Riordan, J. Spencer and G. Tusn\'ady,
\lq\lq The degree sequence of a scale free random graph process," 
{\sl Random Structures and Algorithms} {\bf 18} (2001), pp. 279-290.

\bibitem{CCDM}
G. Caldarelli, A. Capocci, P. De Los Rios and M.A. Mu\~noz,
\lq\lq Scale-free networks from varying vertex intrinsic fitness,"
{\sl Physical Review Letters} {\bf 89} (2002), 258702.

\bibitem{ClausetShaliziNewman} 
A. Clauset, C. Rohilla Shalizi and M.E.J. Newman,
\lq\lq Power-law distributions in empirical data,"
{\sl SIAM Review} {\bf 51} (2009), pp. 661-703.

\bibitem{ChungBook}
K.L. Chung,
{\sl A Course in Probability Theory},
Second Edition, Academic Press,
New York (NY), 1974.

\bibitem{DavidNagarajaBook}
H.A. David and H.N. Nagaraja,
{\sl Order Statistics},
Third Edition,
Wiley Series in Probability and Statistics,
John Wiley \& Sons, Hoboken (NJ),  2003.




\bibitem{Durrett_Book}
R. Durrett,
{\sl Random Graph Dynamics},
Cambridge Series in Statistical and Probabilistic Mathematics,
Cambridge University Press,
Cambridge (UK), 2007.

\bibitem{EKM}
P. Embrechts, C. Kl\"uppelberg and T. Mikosch,
{\sl Modelling Extremal Events for Insurance and Finance}, 
Springer-Verlag, Berlin (Germany), 1997. 

\bibitem{FIKMMU}
A. Fujihara, Y. Ide, N. Konno, N. Masuda, H. Miwa and M. Uchida,
\lq\lq Limit theorems for the average distance  and the degree distribution  of the threshold network model," 
{\em Interdisciplinary Information Sciences} {\bf 15} (2003),
pp. 361-366.

\bibitem{JansonLuczakRucinski}
S. Janson, T. \mbox{\L}uczak and A. Ruci\'{n}ski, 
{\sl Random Graphs},
Wiley-Interscience Series in Discrete Mathematics and Optimization, 
John Wiley \& Sons, New York (NY), 2000.


\bibitem{MakowskiYagan-JSAC}
A. M. Makowski and O. Ya\u{g}an,
\lq\lq Scaling laws for connectivity in random threshold graph models
with non-negative fitness variables,"
{\sl IEEE Journal on Selected Areas in Communications}
{\bf JSAC--31} (2013),
Special Issues on Emerging Technologies in Communications
(Area 4: Social Networks).

\bibitem{NewmanSurvey}
M.E.J. Newman,
\lq\lq The structure and function of complex networks,"
{\sl SIAM Review} {\bf 45} (2003), pp. 167-256.



\bibitem{SPal_Thesis}
S. Pal,
{\sl Adventures on Networks: Degrees and Games},
Ph.D. Thesis, Department of Electrical and Computer Engineering,
University of Maryland, College Park (MD), December 2015.

\bibitem{PalMakowski_CDC2015}
S. Pal and A.M. Makowski,
\lq\lq On the asymptotics of degree distributions,"
in the Proceedings of  the 53rd IEEE Conference on Decision and Control (CDC 2015),
Osaka (Japan), December 2015.

\bibitem{PalMakowski-TNSE}
S. Pal and A.M. Makowski,
\lq\lq Asymptotic distributions in large (homogeneous) random networks: A little theory and a counterexample,"
{\sl IEEE Transactions on Network Science and Engineering}. Accepted for publication, July 2019.
Also available at arXiv:1710.11064.

\bibitem{SC}
V.D.P. Servedio and G. Caldarelli,
\lq\lq Vertex intrinsic fitness: 
How to produce arbitrary scale-free networks,"
{\sl Physical Review E} {\bf 70} (2004), 056126

\bibitem{Shiryayev}
A.N. Shiryayev,
{\sl Probability}, Graduate Texts in Mathematics {\bf 95},
Translated by R.P. Boas, Springer-Verlag, New York (NY), 1984.

\end{thebibliography}
\end{document}